%% file: root.tex
\documentclass[letterpaper, 10pt, final, journal, twocolumn]{IEEEtran}

\input{myheader}

\graphicspath{{./Figures}}
\DeclareGraphicsExtensions{.pdf,.jpg,.png}

\title{Energy-based Modeling of AC Motors}
\author{Pascal Combes,
		Al~Kassem~Jebai,
		François~Malrait,
		Philippe~Martin,
		and~Pierre~Rouchon
	\thanks{P.~Martin and P.~Rouchon are with Centre Automatique et Systèmes, MINES ParisTech, PSL Research University, 75006 Paris, France
	{\tt\footnotesize\{philippe.martin,pierre.rouchon\}@mines-paristech.fr}}
	\thanks{P.~Combes, A.K.~Jebai and F.~Malrait are with Schneider Toshiba Inverter Europe, 27120 Pacy-sur-Eure, France
	{\tt\footnotesize \{pascal.combes, al-kassem. jebai, francois.malrait\}@schneider-electric.com}}}

\begin{document}
	\maketitle
	
	\begin{abstract}
		We propose an approach to modeling of AC motors entirely based on analytical mechanics. Symmetry and connection constraints are moreover incorporated in the energy function from which the models are derived. The approach is especially suited to handle magnetic saturation, but also directly recovers the standard unsaturated models of the literature. The theory is illustrated by some experimental data.
	\end{abstract}
	
	\input{intro}
	\input{mechanics}
	\input{frames}
	\input{symmetries}
	\input{connections}
	\input{simpleModels}

	\input{experiments}

	\section{Conclusion}
		We have proposed an approach to modeling of AC motors entirely based on analytical mechanics, and taken into account geometric symmetries and connection constraints in the energy function. The approach is well-suited to handle magnetic saturation, which is paramount in the design of sensorless control laws are very low speed. In particular, it justifies the fact that star-connected saturated models can be expressed in the $\DQ$ (or $\dq$) frame, which is not obvious though usually taken for granted in the literature.
		
		A further step in modeling would be to incorporate some hysteresis effects, which seem to be the next major phenomenon in AC motors without permanent magnets (SynRM and IM) once magnetic saturation has been taken into account.

	\bibliographystyle{IEEEtran}
	\bibliography{biblio}
	
\end{document}

%% file: myheader.tex
\usepackage[utf8]{inputenc}
\usepackage[T1]{fontenc}
\usepackage[english]{babel}
\usepackage{cite}
\usepackage{times}
\usepackage{url}

\usepackage{graphicx}
\usepackage[caption=false,font=footnotesize]{subfig}
\usepackage{tikz}
\usepackage{pgfplots}
\pgfplotsset{compat=newest}
\pgfplotsset{
	plot coordinates/math parser=false,
	tick label style={%
		/pgf/number format/precision=3,%
		/pgf/number format/std=-3:3},%
	ylabel absolute/.style={%
		every axis y label/.style={at={(0,0.5)},xshift=-25pt,rotate=90},
		every y tick scale label/.style={
			at={(0,1)},above right,inner sep=0pt,yshift=0.3em
		},
	}
}
\newlength\figurewidth
\newlength\figureheight

\usepackage{amsmath}
\usepackage{amssymb}

\usepackage{array}
\usepackage{booktabs}

\usepackage{cite}

\usepackage{siunitx}
\DeclareSIUnit\rpm{rpm}

\usepackage{url}

\usepackage{xcolor}

\def\ordfield{}
\def\subfield{}
\def\expfield{}
\def\getfields#1{%
	\def\ordfield{}%
	\def\subfield{}%
	\def\expfield{}%
	\expandafter\getexpfield\expandafter\beginsym#1^!\endsymb%
}
\def\extractsubfield\beginsymb#1_!\endsymb{%
	\def\subfield{#1}%
}
\def\getsubfield\beginsymb#1_#2\endsymb{%
	\if!#2\relax\else%
	\extractsubfield\beginsymb#2\endsymb%
	\fi%
	\def\ordfield{#1}%
}
\def\extractexpfield\beginsymb#1^!\endsymb{%
	\getsubfield\beginsymb#1_!\endsymb%
	\expandafter\def\expandafter\expfield\expandafter{\ordfield}
}
\def\getexpfield\beginsym#1^#2\endsymb{%
	\if!#2\relax\else%
	\extractexpfield\beginsymb#2\endsymb%
	\fi%
	\getsubfield\beginsymb#1_!\endsymb%
}

\newcommand{\be}{\begin{equation}}
\newcommand{\ee}{\end{equation}}
\newcommand{\beno}{\begin{equation*}}
\newcommand{\eeno}{\end{equation*}}
\newcommand{\ba}{\begin{IEEEeqnarray}{rCl}\IEEEyesnumber}
\newcommand{\ea}{\end{IEEEeqnarray}}
\newcommand{\bano}{\begin{IEEEeqnarray*}{rCl}}
\newcommand{\eano}{\end{IEEEeqnarray*}}
\newcommand{\ans}{\IEEEeqnarraynumspace}
\newcommand{\ase}{\IEEEyessubnumber}
\newcommand{\ane}{\IEEEnonumber}

\newcommand{\sectionname}{section}
\def\firstuppercase#1{\uppercase{#1}}

\newcommand{\eqnref}[1]{\eqref{eqn:#1}}
\newcommand{\secref}[1]{\sectionname~\ref{sec:#1}}
\newcommand{\Secref}[1]{\expandafter\firstuppercase\sectionname~\ref{sec:#1}}
\newcommand{\figref}[1]{Fig.~\ref{fig:#1}}
\newcommand{\Figref}[1]{\expandafter\firstuppercase\figurename~\ref{fig:#1}}
\newcommand{\tblref}[1]{\expandafter\lowercase\expandafter{\tablename}~\ref{tbl:#1}}

\makeatletter
\newcommand{\transpose}[1]{#1^T}
\newcommand{\Transpose}[1]{{#1}^T}

\newcommand{\tderive}[2][1]{%
	\def\tempa{#1}%
	\def\tempb{1}%
	\def\tempc{2}%
	\ifx\tempb\tempa%
		\getfields{#2}%
		\dot{\ordfield}%
		\ifx\expfield\empty\else%
			^{\expfield}%
		\fi%
		\ifx\subfield\empty\else%
			_{\subfield}%
		\fi
	\else%
		\ifx\tempc\tempa%
			\getfields{#2}%
			\ddot{\ordfield}%
			\ifx\expfield\empty\else%
				^{\expfield}%
			\fi%
			\ifx\subfield\empty\else%
				_{\subfield}%
			\fi
		\else%
			#2^{(#1)}%
		\fi%
	\fi%
}
\newcommand{\Tderive}[2][1]{%
	\def\tempa{#1}%
	\def\tempb{1}%
	\ifx\tempa\tempb%
		\frac{d#2}{dt}%
	\else%
		\frac{d^{#1}#2}{dt^{#1}}%
	\fi%
}
\newcommand{\pderive}[3][1]{%
	\def\tempa{#1}%
	\def\tempb{1}%
	\ifx\tempa\tempb%
		\partial_{#3}#2%
	\else%
		\partial^{#1}_{#3}#2%
	\fi%
}
\newcommand{\Pderive}[2]{%
	\nabla_{#2}#1%
}
\newcommand{\Jac}[2]{%
	\partial_{#2}#1
}
\newcommand{\PPderive}[3][!]{%
	\if!#1%
		\frac{\partial^2#2}{\partial{#3}^2}%
	\else%
		\frac{\partial^2#1}{\partial{#2}\partial{#3}}%
	\fi%
}
\newcommand{\ind}[3]{%
	\def\tempa{#3}%
	\ifx\tempa\empty%
		#1_{#2}%
	\else%
		#1_{#2#3}%
	\fi%
}
\newcommand{\fcos}[1]{a_0^{#1}}
\newcommand{\fsin}[1]{b_0^{#1}}
\newcommand{\fscos}[1]{a_{0\star}^{#1}}
\newcommand{\fssin}[1]{b_{0\star}^{#1}}

\newcommand{\lf}[1]{%
	\getfields{#1}%
	\bar{\ordfield}%
	\ifx\expfield\empty\else%
		^{\expfield}%
	\fi%
	\ifx\subfield\empty\else%
		_{\subfield}%
	\fi
}
\newcommand{\hf}[1]{
	\getfields{#1}%
	\widetilde{\ordfield}%
	\ifx\expfield\empty\else%
		^{\expfield}%
	\fi%
	\ifx\subfield\empty\else%
		_{\subfield}%
	\fi
}
\makeatother

\newcommand{\abc}{abc}
\newcommand{\ABZ}{\alpha\beta0}
\newcommand{\AB}{\alpha\beta}
\newcommand{\dqZ}{dq0}
\newcommand{\dq}{dq}
\newcommand{\DQZ}{DQ0}
\newcommand{\DQ}{DQ}
\newcommand{\xyz}{xyz}
\newcommand{\xyZ}{xy0}
\newcommand{\xy}{xy}

\newcommand{\Lg}{\mathcal{L}}
\newcommand{\Labc}{\Lg^{\abc}}
\newcommand{\Hm}{\mathcal{H}}
\newcommand{\Habc}{\Hm^{\abc}}
\newcommand{\HABZ}{\Hm^{\ABZ}}
\newcommand{\HdqZ}{\Hm^{\dqZ}}
\newcommand{\HmdqZ}{\Hm_m^{\dqZ}}
\newcommand{\HDQZ}{\Hm^{\DQZ}}
\newcommand{\HxyZ}{\Hm^{\xyZ}}
\newcommand{\Hxyz}{\Hm^{\xyz}}
\newcommand{\HmDQZ}{\Hm_m^{\DQZ}}

\newcommand{\HcxyZ}{\Hm_c^{\xyZ}}
\newcommand{\Hcsxy}{\Hm_{c\star}^{\xy}}

\newcommand{\Hzxyz}{\Hm_0^{\xyz}}

\newcommand{\Hzsxy}{\Hm_{0\star}^{\xy}}
\newcommand{\HmzDQZ}{\Hm_{m0}^{\DQZ}}
\newcommand{\HmzsDQ}{\Hm_{m0\star}^{\DQ}}
\newcommand{\HsxyZ}{\Hm_\star^{\xyZ}}

\newcommand{\HbmZDQZ}{\overline{\Hm}^{\DQZ}_{m0}}
\newcommand{\HbmZsDQ}{\overline{\Hm}^{\DQ}_{m0\star}}

\newcommand{\HDQ}{\Hm^{\DQ}}

\newcommand{\zero}[2]{0_{#1,#2}}
\newcommand{\Id}{I_}
\newcommand{\JJ}{\mathcal{J}_}
\newcommand{\Rot}{\mathcal{R}_}
\newcommand{\CC}{\mathcal{C}}
\newcommand{\Pabc}{\mathcal{P}^{\abc}}

\newcommand{\Oabc}{\mathcal{O}^{\abc}}
\newcommand{\OABZ}{\mathcal{O}^{\ABZ}}

\newcommand{\qsabc}{q_s^{\abc}}
\newcommand{\qrabc}{q_r^{\abc}}

\newcommand{\is}{\imath_s}
\newcommand{\ir}{\imath_r}
\newcommand{\isabc}{\is^{\abc}}
\newcommand{\irabc}{\ir^{\abc}}
\newcommand{\isa}{\is^a}
\newcommand{\isb}{\is^b}
\newcommand{\isc}{\is^c}
\newcommand{\ira}{\ir^a}
\newcommand{\irb}{\ir^b}
\newcommand{\irc}{\ir^c}
\newcommand{\isABZ}{\is^{\ABZ}}
\newcommand{\isAB}{\is^{\AB}}
\newcommand{\irABZ}{\ir^{\ABZ}}
\newcommand{\isA}{\is^\alpha}
\newcommand{\irA}{\ir^\alpha}
\newcommand{\isB}{\is^\beta}
\newcommand{\irB}{\ir^\beta}
\newcommand{\isZ}{\is^0}
\newcommand{\irZ}{\ir^0}
\newcommand{\isdqZ}{\is^{\dqZ}}
\newcommand{\irdqZ}{\ir^{\dqZ}}
\newcommand{\isdq}{\is^{\dq}}
\newcommand{\irdq}{\ir^{\dq}}
\newcommand{\isd}{\is^d}
\newcommand{\isq}{\is^q}
\newcommand{\ird}{\ir^d}
\newcommand{\irq}{\ir^q}
\newcommand{\isDQZ}{\is^{\DQZ}}
\newcommand{\irDQZ}{\ir^{\DQZ}}
\newcommand{\isDQ}{\is^{\DQ}}

\newcommand{\isD}{\is^D}
\newcommand{\isQ}{\is^Q}
\newcommand{\irD}{\ir^D}
\newcommand{\irQ}{\ir^Q}
\newcommand{\isxyZ}{\is^{\xyZ}}
\newcommand{\irxyZ}{\ir^{\xyZ}}
\newcommand{\isxy}{\is^{\xy}}
\newcommand{\irxy}{\ir^{\xy}}
\newcommand{\isxyz}{\is^{\xyz}}
\newcommand{\irxyz}{\ir^{\xyz}}

\newcommand{\phis}{\lambda_s}
\newcommand{\phir}{\lambda_r}
\newcommand{\phisabc}{\phis^{\abc}}
\newcommand{\phirabc}{\phir^{\abc}}
\newcommand{\phisa}{\phis^a}
\newcommand{\phisb}{\phis^b}
\newcommand{\phisc}{\phis^c}

\newcommand{\phisABZ}{\phis^{\ABZ}}
\newcommand{\phirABZ}{\phir^{\ABZ}}
\newcommand{\phisA}{\phis^\alpha}
\newcommand{\phirA}{\phir^\alpha}
\newcommand{\phisB}{\phis^\beta}
\newcommand{\phirB}{\phir^\beta}
\newcommand{\phisZ}{\phis^0}
\newcommand{\phirZ}{\phir^0}
\newcommand{\phisdqZ}{\phis^{\dqZ}}
\newcommand{\phirdqZ}{\phir^{\dqZ}}
\newcommand{\phisdq}{\phis^{\dq}}
\newcommand{\phirdq}{\phir^{\dq}}
\newcommand{\phisd}{\phis^d}
\newcommand{\phisq}{\phis^q}
\newcommand{\phird}{\phir^d}
\newcommand{\phirq}{\phir^q}
\newcommand{\phisDQZ}{\phis^{\DQZ}}
\newcommand{\phirDQZ}{\phir^{\DQZ}}
\newcommand{\phisDQ}{\phis^{\DQ}}

\newcommand{\phisD}{{\phis^D}}
\newcommand{\psisD}{{\psi_s^D}}
\newcommand{\phisQ}{{\phis^Q}}
\newcommand{\phirD}{\phir^D}
\newcommand{\phirQ}{\phir^Q}
\newcommand{\phisxyZ}{\phis^{\xyZ}}
\newcommand{\phirxyZ}{\phir^{\xyZ}}
\newcommand{\phisxy}{\phis^{\xy}}
\newcommand{\phirxy}{\phir^{\xy}}
\newcommand{\phisxyz}{\phis^{\xyz}}
\newcommand{\phirxyz}{\phir^{\xyz}}

\newcommand{\us}{u_s}
\newcommand{\usabc}{\us^{\abc}}
\newcommand{\usa}{\us^a}
\newcommand{\usb}{\us^b}
\newcommand{\usc}{\us^c}
\newcommand{\usABZ}{\us^{\ABZ}}
\newcommand{\usxyz}{\us^{\xyz}}
\newcommand{\usxyZ}{\us^{\xy0}}
\newcommand{\usxy}{\us^{\xy}}
\newcommand{\usZ}{\us^0}
\newcommand{\usdqZ}{\us^{\dqZ}}
\newcommand{\usdq}{\us^{\dq}}
\newcommand{\usDQZ}{\us^{\DQZ}}
\newcommand{\usDQ}{\us^{\DQ}}

\newcommand{\vs}{v_s}
\newcommand{\vsabc}{\vs^{\abc}}
\newcommand{\vsa}{\vs^a}
\newcommand{\vsb}{\vs^b}
\newcommand{\vsc}{\vs^c}
\newcommand{\vsxy}{\vs^{\xy}}
\newcommand{\vsDQ}{\vs^{\DQ}}
\newcommand{\vsAB}{\vs^{\AB}}
\newcommand{\vsZ}{\vs^0}

\newcommand{\te}{\theta}
\newcommand{\tm}{\te_m}
\newcommand{\ts}{\te_s}
\newcommand{\we}{\omega}

\newcommand{\ws}{\we_s}
\newcommand{\ke}{\rho}

\newcommand{\Te}{T_e}
\newcommand{\Teabc}{T_e^{\abc}}
\newcommand{\TeABZ}{T_e^{\ABZ}}
\newcommand{\TeAB}{T_e^{\AB}}
\newcommand{\TedqZ}{T_e^{\dqZ}}
\newcommand{\TeDQZ}{T_e^{\DQZ}}
\newcommand{\Texyz}{T_e^{\xyz}}
\newcommand{\TexyZ}{T_e^{\xy0}}
\newcommand{\Tedq}{T_e^{\dq}}
\newcommand{\TeDQ}{T_e^{\DQ}}
\newcommand{\Texy}{T_e^{\xy}}
\newcommand{\Tl}{T_l}

\newcommand{\Rs}{R_s}
\newcommand{\Rr}{R_r}

\newcommand{\Gs}{\Gamma_{\!s}}
\newcommand{\Gm}{\Gamma_{\!m}}
\newcommand{\Gls}{\Gamma_{\!\!ls}}
\newcommand{\Glr}{\Gamma_{\!\!lr}}

\newcommand{\GsD}{\Gs^D}
\newcommand{\GsQ}{\Gs^Q}
\newcommand{\GsZ}{\Gs^0}

\newcommand{\GlsZ}{\Gls^0}
\newcommand{\GlrZ}{\Glr^0}

\newcommand{\phiM}{\Phi_M}

\newcommand{\starpoint}{N}

\newcommand{\fD}{{\rm f}^D}
\newcommand{\fQ}{{\rm f}^Q}
\newcommand{\fx}{{\rm f}^X}

\newcommand{\Jl}{J}
\newcommand{\np}{n}
\newcommand{\nr}{n_r}

\newif\ifspmyesipmno
\spmyesipmnofalse


%% file: intro.tex
\section{Introduction}\label{sec:intro}
\IEEEPARstart{G}{ood} models are usually paramount to design good control laws. This is the case for AC motors ---Permanent Magnet Synchronous Motors (PMSM), Synchronous Reluctance Motors (SynRM), Induction Motors (IM),~\ldots---, especially when ``sensorless'' control is considered. In this mode of operation, neither the rotor position nor its velocity is measured, and the control law must make do with only current measurements; a suitable model is therefore essential to relate the currents to the other variables. When operating above moderately low speed, i.e., above about $10\%$ of the rated speed, models neglecting magnetic saturation are usually accurate enough for control purposes; but at low speed, magnetic (cross-)saturation must be taken into account, in particular when high-frequency signal injection is used, see e.g.~\cite{LiZHBS2009ITIA,SergeDM2009ITM,BiancFB2011ECCE,JebaiMMR2016IJoC}.

The usual approach is to first derive an unsaturated sinusoidal model of the motor by computing the magnetomotive force and considering linear flux-current relations, see e.g.~\cite{Vas1998book,KrauseWSP2013book,Chiasson2005book}. If necessary, the model is then extended to take into account magnetic saturation by using nonlinear flux-current relations, see e.g.~\cite{Vas1998book}, and nonsinusoidal windings by adding harmonics to the back-emf and torque, see e.g.~\cite{HoltzS1996ITIE,PetroOST2000ITPE}. Besides being rather tedious to carry out ---even for the unsaturated case---, this approach may be awkward to handle saturation. Indeed the flux-current relations should respect the so-called reciprocity conditions~\cite{MelkeW1990IAITo}; this is always the case for linear flux-current relations, but can be troublesome to enforce in more involved cases. Moreover, it is not easy to correctly describe the impact of saturation on the electromagnetic torque. Finally, it usually requires a good knowledge of the motor internal layout.

An alternative route is to use analytical mechanics, which can be seen as a ``top-down'' approach compared to the usual ''bottom-up" approach. Analytical mechanics models any physical system by a single scalar function, related to energy, which conveys all the required information to obtain a model. It does not bring any new information with respect to first-order principles (d'Alembert principle, Kirchoff laws,~\ldots), but is often a convenient formulation when dealing with complex systems. Energy-based modeling is not a new method, neither in electromechanical systems~\cite{WhiteW1959book}, nor in control theory~\cite{JeltsS2009ICSM}. For AC motors, it is even customary to compute the electromagnetic torque by some energy/coenergy considerations, see e.g.~\cite{Vas1998book,KrauseWSP2013book,Chiasson2005book}, though it remains an aid rather than a full application of analytical mechanics.

In this paper, which extends the preliminary work~\cite{JCMMR2014CDC}, we derive models of AC motors by relying solely on analytical mechanics. The current-flux relations as well as the electromagnetic torque expression are expressed as partial derivatives of a single ``energy function'' which encodes all the required information, such as magnetic saturation. Moreover, we take advantage of basic geometric symmetries enjoyed by any well-built motor, and of constraints induced by the connections of the motor windings, to simplify the models. This provides an elegant and efficient way to model AC motors, with several benefits:
\begin{itemize}
	\item  the unsaturated sinusoidal models are easily retrieved as deriving from the simplest possible energy functions, without requiring a precise knowledge of the motor internal layout
	\item more importantly, it produces models which are by construction ``physically correct'', even in the saturated non-sinusoidal case: the flux-current relations respect the reciprocity conditions and are coherent with the electromagnetic torque, since all these quantities are derived from a single energy function
	\item last but not least, it justifies the fact that a (star-connected) saturated model can be expressed in the $\DQ$ (or $\dq$) frame. Indeed, saturation physically takes place in the $\abc$ frame, and it is not obvious that the corresponding nonlinearities behave well under a change of frames and that the $0$-axis can be decoupled. Though not completely obvious, this issue is very seldom addressed in the literature.
\end{itemize}
The approach turns out to be very effective in practice to model saturated motors, since all the effort consists in suitably ``shaping'' a single scalar function. This is illustrated and experimentally validated in \secref{PMSM} on an actual PMSM.

The paper runs as follows: in \secref{meca}, we recall the basics of analytical mechanics, and apply them to AC motors; in \secref{frames}, we show how the useful standard frame changes operate in our formalism; in \secref{symmetries}, we show how the various geometric symmetries of AC motors constrain, hence simplify, their energy functions; in \secref{connections}, we do the same with the connection constraints of the stator and rotor windings; in \secref{models}, we then put together all the material of the previous sections to effortlessly recover the classical unsaturated sinusoidal models; finally, we treat in~\secref{PMSM} the case of a real PMSM with non-sinusoidal and saturation effects, including comparisons with experimental data.

We end this section with a few notations that will be used throughout the paper. The coordinates of the $3\times1$ vector $x$ in the $uvw$ frame are denoted by $x^{uvw} := \transpose{\left(x^u, x^v, x^w\right)}$. The Jacobian 
matrix of the (scalar or vector) function~$f$ is $\Jac{f}{x^{uvw}}:=\left(\Jac{f}{x^u}, \Jac{f}{x^v}, \Jac{f}{x^w}\right)$, and its gradient is $\Pderive{f}{x^{uvw}} := \transpose{(\Jac{f}{x^{uvw}})}$; the Hessian of the scalar function~$g$ is thus $\Jac{\Pderive{g}{x^{uvw}}}{x^{uvw}}$. We will also use the following matrices:  $\Id n$, the $n\times n$~identity matrix, and
\beno
\JJ2 := \begin{pmatrix}
	0 & -1 \\
	1 & 0
\end{pmatrix} \quad 
\Rot2(\eta) := \begin{pmatrix}
	\cos{\eta} & -\sin{\eta} \\
	\sin{\eta} & \cos{\eta}
\end{pmatrix}	
\eeno
\beno
\JJ3 := \begin{pmatrix}
	0 & -1 & 0 \\
	1 & 0 & 0 \\
	0 & 0 & 0
\end{pmatrix} 
\quad	\Rot3(\eta) := \begin{pmatrix}
	\cos{\eta} & -\sin{\eta} & 0 \\
	\sin{\eta} & \cos{\eta} & 0 \\
	0 & 0 & 1
\end{pmatrix};
\eeno
notice the useful relation $\frac{d\Rot i}{d\eta}(\eta)=\JJ i\Rot i(\eta)$ for $i=2,3$.

%% file: mechanics.tex
\section{Analytical mechanics and AC motors}\label{sec:meca}

\subsection{Basics of analytical mechanics (see e.g.~\cite{Lurie2002book})}\label{sec:basics}
In analytical mechanics, the time evolution of a physical system (mechanical, electrical, electromechanical, \dots) is described by the Euler-Lagrange equations,
\ba\label{eqn:meca:dynamic}
	\Tderive{q} &=& \tderive{q} \ase \label{eqn:meca:dynamic:qdot} \\
	\Tderive{} \Pderive{\Lg}{\tderive{q}}(q, \tderive{q}) &=& \Pderive{\Lg}{q}(q, \tderive{q}) + \mathcal{F}; \ase\label{eqn:meca:dynamic:q}
\ea
the components of~$q\in\mathbb{R}^n$ are called the generalized coordinates, and those of~$\tderive{q}\in\mathbb{R}^n$ the generalized velocities. The scalar function~$\Lg$, called the Lagrangian of the system, encodes all the necessary information; in full generality it may also depend on the time~$t$, though it is not needed in this paper. Finally, the vector function~$\mathcal{F}$ represents the generalized forces acting on the system. 

One drawback of the Lagrangian formulation~\eqnref{meca:dynamic} is that it does not yield a system of equations in state form, which is not very convenient for control purposes. When $\Lg$~is not degenerate, i.e., $\pderive{}{\tderive{q}}\Pderive{\Lg}{\tderive{q}}$~is full rank, the state form associated to~\eqnref{meca:dynamic} can be obtained by the following Hamiltonian formulation. The new scalar function~$\Hm$, called the Hamiltonian of the system, is defined by
\be
	\Hm(p, q) := \transpose{\tderive{q}}p - \Lg(q, \tderive{q}), \label{eqn:meca:legendre}
\ee
where $p := \Pderive{\Lg}{\tderive{q}}(q, \tderive{q}) \in \mathbb{R}^n$~is the vector of so-called generalized momenta.
In \eqnref{meca:legendre}, $\tderive{q}$~must be expressed in terms of~$p$ and~$q$ by inverting the relation~$p = \Pderive{\Lg}{\tderive{q}}(q, \tderive{q})$. It is not necessary to perform this inversion explicitly to find the state form. Indeed, we find using~\eqnref{meca:legendre} that the differential of~$\Hm$ is 
\bano
	d\Hm &=& \transpose{p}d\tderive{q} + \transpose{\tderive{q}}dp - \Jac{\Lg}{q}dq - \Jac{\Lg}{\tderive{q}}d\tderive{q} \ane \\ 
		&=& \transpose{\tderive{q}}dp - \Jac{\Lg}{q}dq.
\eano
Since by definition 
\beno
	d\Hm = \Jac{\Hm}{p}dp + \Jac{\Hm}{q}dq,
\eeno
this immediately gives the state equations
\ba\label{eqn:meca:hamiltonian}
	\Tderive{p} &=& -\Pderive{\Hm}{q}(p, q) + \mathcal{F} \ase \label{eqn:meca:hamiltonian:p} \\
	\Tderive{q} &=& \Pderive{\Hm}{p}(p, q), \ase \label{eqn:meca:hamiltonian:q}
\ea
by identifying the two expressions for~$d\Hm$ and using~\eqref{eqn:meca:dynamic}. Notice in particular that $\tderive{q} = \Pderive{\Hm}{p}(p, q)$, which gives the inverse of the relation $p = \Pderive{\Lg}{\tderive{q}}(q, \tderive{q})$. Moreover, differentiating \
\beno
	\tderive{q} = \Pderive{\Hm}{p}\bigl(\pderive{\Lg}{\tderive{q}}(q, \tderive{q}), q\bigr)
\eeno
with respect to~$\tderive{q}$ immediately yields
\beno
	\pderive{}{p}\Pderive{\Hm}{p} \cdot \pderive{}{\tderive{q}}\Pderive{\Lg}{\tderive{q}} = \Id{};
\eeno
hence the nondegeneracy condition~$\pderive{}{\tderive{q}}\Pderive{\Lg}{\tderive{q}}$ full rank is equivalent to~$\pderive{}{p}\Pderive{\Hm}{p}$ full rank.

Both the Lagrangian and the Hamiltonian have the dimension of energy. However, since
\ba
	\Tderive{\Lg} &=& \Jac{\Lg}{q}\Tderive{q} + \Jac{\Lg}{\tderive{q}}\Tderive{\tderive{q}} \ane \\ 
		&=& \Jac{\Lg}{q}\tderive{q} + \Tderive{} \left(\Jac{\Lg}{\tderive{q}}\tderive{q}\right) - \Tderive{} \left(\Jac{\Lg}{\tderive{q}}\right)\tderive{q} \ane \\ 
		&=& \Tderive{} \left(\Jac{\Lg}{\tderive{q}}\tderive{q}\right) - \transpose{\mathcal{F}}\tderive{q} \ane \label{eqn:meca:dL} \\ 
	\Tderive{\Hm} &=& \Tderive{} (\transpose{\tderive{q}}p) - \Tderive{\Lg} \ane \\
		&=& \transpose{\mathcal{F}}\tderive{q}, \ane \label{eqn:meca:d}
\ea
the Hamiltonian is truly the energy of the system, whereas the Lagangian has no physical meaning. Indeed the variation of~$\Hm$ over time is~$\transpose{\mathcal{F}}\tderive{q}$, the work of the generalized forces along the trajectory of the system; in particular, $\Hm$ is conserved when $\mathcal{F} = 0$.

\subsection{Application to AC motors}\label{sec:appli}
The formulation presented in the previous section is now applied to a general three-phase AC motor. A motor with $\np$~pairs of poles has 3$\np$~identical stator windings, and we classically assume that the rotor can be modeled by 3$\np$~identical windings as well (this idealization is valid in particular for squirrel-cage IMs).
As all the poles are built identically, the motor can be reduced to a single pair of poles by considering the electrical rotor angle~$\te := \np\tm$ instead of the mechanical rotor angle~$\tm$. The three windings of the equivalent simplified motor are labeled~$a$, $b$~and~$c$.
For such an electromechanical system, the generalized coordinates are
\beno
q := (\qsabc, \qrabc, \te),
\eeno
where $\qsabc$~is the vector of charges in the frame defined by the stator windings, and $\qrabc$,~the vector of charges in the frame defined by the rotor windings. These coordinates are independent, as we consider for the time being that the windings are unconnected; the connection constraints will be dealt with in \secref{connections}, see Figs.~\ref{fig:connection:rotor}~and~\ref{fig:connection:stator}. The corresponding generalized velocities are
\beno
\tderive{q} := (\isabc, \irabc, \we),
\eeno
where in view of~\eqnref{meca:dynamic:qdot} $\isabc$ is the vector of currents in the stator windings, $\irabc$ is the vector of currents in the rotor windings, and $\we$ is the (electrical) angular velocity.
As there is no storage of electrical charges in the windings, the Lagrangian of the motor does not depend on the charges, i.e., has the form
\beno
\Labc(\isabc, \irabc, \te, \we). 
\eeno

\begin{figure}
	\centering
	\includegraphics[width=250pt]{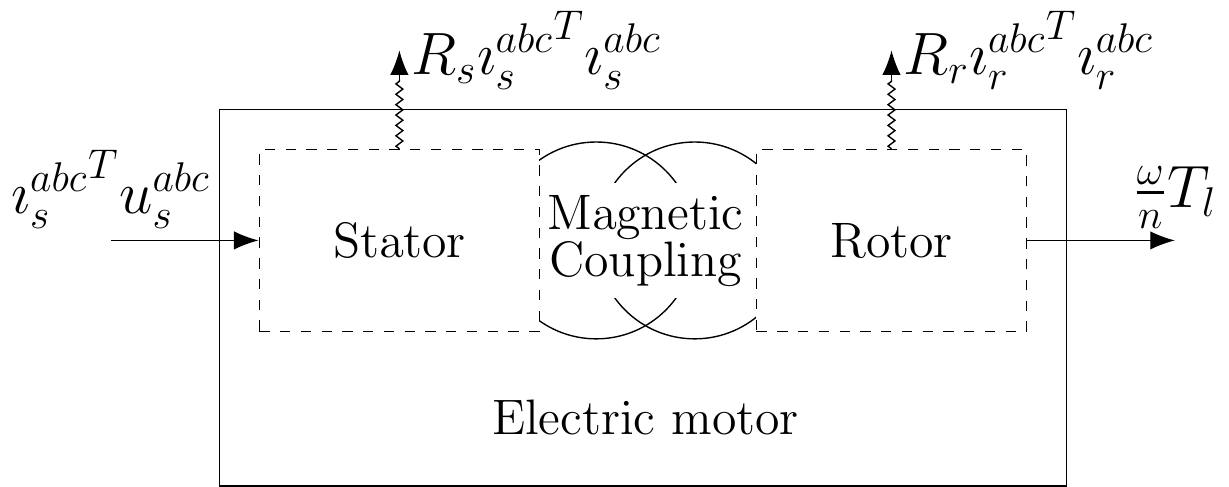}
	\caption{Power exchanges in an electric motor.}
	\label{fig:system}
\end{figure}

To compute the generalized forces, we notice that according to \figref{system} the motor:
\begin{itemize}
	\item receives the electrical power~$\Transpose{\isabc}\usabc$, where $\usabc$~is the vector of voltage drops across the stator windings; hence, $\usabc$ is a generalized force associated to the generalized coordinate~$\qsabc$
	\item loses the electrical power~$\Rs\Transpose{\isabc}\isabc$ in its stator resistance~$\Rs$; hence, $-\Rs\isabc$~is a generalized force associated to the generalized coordinate~$\qsabc$
	\item loses the electrical power~$\Rr\Transpose{\irabc}{\irabc}$ in the rotor resistance~$\Rr$; hence, $-\Rr\irabc$~is the generalized force associated to the generalized coordinate~$\qrabc$
	\item yields the mechanical power~$\frac{\we}{\np}\Tl$, where $\Tl$~is the load torque; hence, $-\frac{\Tl}{\np}$~is the generalized force associated to the generalized coordinate~$\te$.
\end{itemize}
Taking into account these generalized forces, the Euler-Lagrange equations~\eqnref{meca:dynamic} read
\ba\label{eqn:appli:euler-lagrange}
\Tderive{\te} &=& \we \ase \\
\Tderive{} \Pderive{\Labc}{\isabc}(\isabc, \irabc, \te, \we) &=& \usabc - \Rs\isabc \ase \\
\Tderive{} \Pderive{\Labc}{\irabc}(\isabc, \irabc, \te, \we) &=& -\Rr\irabc \ase \\
\Tderive{} \Pderive{\Labc}{\we}(\isabc, \irabc, \te, \we) &=& \frac{\Te}{\np} - \frac{\Tl}{\np}, \ase 
\ea
where the electromagnetic torque~$\Te$ is defined by
\be
\Te := \np\Pderive{\Labc}{\te}(\isabc, \irabc, \te, \we); \label{eqn:appli:lagrangian:te}
\ee
as~$\Labc$ does not depend on~$\qsabc$ and~$\qrabc$, it is useless to include in~\eqnref{appli:euler-lagrange} the two equations
\bano
\Tderive{\qsabc} &=& \isabc \\
\Tderive{\qrabc} &=&\irabc.
\eano
The system \eqref{eqn:appli:euler-lagrange}-\eqref{eqn:appli:lagrangian:te} is then completely defined once the scalar function~$\Labc$ is specified.

To get state equations we define the generalized momenta $p:=(\phisabc, \phirabc, \ke)$ associated to the generalized velocities by
\ba
\phisabc &:=& \Pderive{\Labc}{\isabc}(\isabc, \irabc, \te, \we) \ase \label{eqn:appli:iphi:stator} \\
\phirabc &:=& \Pderive{\Labc}{\irabc}(\isabc, \irabc, \te, \we) \ase \label{eqn:appli:iphi:rotor} \\
\ke &:=& \np^2 \Pderive{\Labc}{\we}(\isabc, \irabc, \te, \we); \ans \ase \label{eqn:appli:lagrangian:w}
\ea
$\phisabc$~and $\phirabc$~are the flux linkages of the stator and rotor windings, and $\ke$~is the kinetic momentum.
The $\np^2$~factor in \eqnref{appli:lagrangian:w} stems from the fact that the transformation~$\te = \np\tm$ is not normalized, so
\beno
\ke = \np \Pderive{\Labc}{\tderive{\tm}} = \np \Pderive{\Labc}{\we}\cdot\Pderive{\we}{\tderive{\tm}} = \np^2 \Pderive{\Labc}{\we}.
\eeno
As in~\eqnref{meca:legendre}, we define the Hamiltonian $\Habc$ by
\beno
\Habc(\phisabc, \phirabc, \te, \ke) := \Transpose{\isabc}\!\phisabc + \we\ke - \Labc(\isabc, \irabc, \te, \ke).
\eeno
The state equations~\eqref{eqn:meca:hamiltonian} then read
\ba\label{eqn:appli:state-form}
\Tderive{\phisabc} &=& \usabc - \Rs\isabc \ase  \label{eqn:appli:state-form:phis} \\
\Tderive{\phirabc} &=& - \Rr\irabc \ase \\
\Tderive{\te} &=& \we \ase \\
\frac{1}{\np} \Tderive{\ke} &=& \Te - \Tl, \ase
\ea
with
\ba\label{eqn:appli:algebraic}
\isabc &=& \Pderive{\Habc}{\phisabc}(\phisabc, \phirabc, \te, \ke) \ase \label{eqn:appli:phii:stator} \\
\irabc &=& \Pderive{\Habc}{\phirabc}(\phisabc, \phirabc, \te, \ke) \ase \label{eqn:appli:phii:rotor} \\
\we &=& \np^2\Pderive{\Habc}{\ke}(\phisabc, \phirabc, \te, \ke) \ans \ase \label{eqn:appli:hamiltonian:w} \\
\Te &=& -\np\Pderive{\Habc}{\te}(\phisabc, \phirabc, \te, \ke); \ans \ase \label{eqn:appli:hamiltonian:te}
\ea
as~$\Habc$ does not depend on~$\qsabc$ and~$\qrabc$, it is useless to include in~\eqnref{appli:state-form} the two equations
\bano
\Tderive{\qsabc} &=& \Pderive{\Habc}{\phisabc}(\phisabc, \phirabc, \te, \ke) \\
\Tderive{\qrabc} &=& \Pderive{\Habc}{\phirabc}(\phisabc, \phirabc, \te, \ke).
\eano
The system \eqref{eqn:appli:state-form}-\eqref{eqn:appli:algebraic} is then completely defined as soon as the scalar function~$\Habc$ is given. 

Notice that $\phisabc$, $\phirabc$, $\ke$ and~$\te$ are the ``natural'' variables to write the system in state form, and turn out to be very convenient for the analysis of models with magnetic saturation, see e.g. section~\ref{sec:nonlin}. It is of course possible to rewrite~\eqref{eqn:appli:euler-lagrange}-\eqref{eqn:appli:lagrangian:te} in state form with the variables~$\isabc$, $\irabc$, $\te$ and~$\we$, but the resulting system is more difficult to manipulate.

%% file: frames.tex
\section{Changing frames}\label{sec:frames}

The physical $abc$~frame is not very convenient for computations. It is therefore customary to rewrite the motor equations in an orthogonal frame. The constant $\abc$-to-$\ABZ$ transformation is first applied. A time-varying transformation is then applied to express the stator and rotor variables in a common frame; though the general case is no more difficult to handle, we consider only the two most widely used transformations:
\begin{itemize}
	\item the $\ABZ$-to-$\dqZ$ transformation, which expresses the stator and rotor variables in the frame rotating at the angular velocity~$\ws=\Tderive{\ts}$, where $\ts$ is the orientation of the control voltage vector; $\ws$~can be seen as a new control. This transformation is commonly used for asynchronous motors 
	\item the $\ABZ$-to-$\DQZ$ transformation, which expresses the stator variables in the rotor frame. This transformation is commonly used for synchronous motors.
\end{itemize}

\subsection{$\abc$-to-$\ABZ$~transformation}\label{sec:ortho}
Define the $\ABZ$~state variables 
\bano
\phisABZ &:=& \CC \phisabc \\
\phirABZ &:=& \CC \phirabc
\eano
and the scalar function
\beno
\HABZ(\phisABZ, \phirABZ, \te, \ke) := \Habc(\transpose{\CC}\phisABZ, \transpose{\CC}\phirABZ, \te, \ke),
\eeno
where $\CC$~is the orthogonal matrix
\beno
\CC := \sqrt{\frac{2}{3}}\begin{pmatrix}
	1 & -\frac{1}{2} & -\frac{1}{2} \\
	0 & \frac{\sqrt{3}}{2} & -\frac{\sqrt{3}}{2} \\
	\frac{\sqrt{2}}{2} & \frac{\sqrt{2}}{2} & \frac{\sqrt{2}}{2}
\end{pmatrix}.
\eeno
We then find
\bano
\Tderive{\phisABZ} &=& \CC \Tderive{\phisabc} \\
&=& \CC (\usabc - \Rs\isabc) \\
&=& \CC \usabc - \Rs\CC\Pderive{\Habc}{\phisabc} \\
&=& \CC \usabc - \Rs\CC\transpose{\CC}\Pderive{\HABZ}{\phisABZ} \\
\Tderive{\phirABZ} &=& - \Rr\CC\transpose{\CC}\Pderive{\HABZ}{\phirABZ}.
\eano
This yields
\ba\label{eqn:ortho:state-form}
\Tderive{\phisABZ} &=& \usABZ - \Rs\isABZ \ase \\
\Tderive{\phirABZ} &=& - \Rr\irABZ \ase \\
\Tderive{\te} &=& \we \ase \\
\frac{1}{\np} \Tderive{\ke} &=& \Te - \Tl \ase,
\ea
where~$\usABZ := \CC\usabc$ and
\ba\label{eqn:ortho:algebraic}
\isABZ &:=& \CC\isabc = \Pderive{\HABZ}{\phisABZ}(\phisABZ, \phirABZ, \te, \ke) \ase \label{eqn:ortho:stator} \\
\irABZ &:=& \CC\irabc = \Pderive{\HABZ}{\phirABZ}(\phisABZ, \phirABZ, \te, \ke) \ase \label{eqn:ortho:rotor} \\
\we &=& \np^2\Pderive{\HABZ}{\ke}(\phisABZ, \phirABZ, \te, \ke) \ase \label{eqn:ortho:w} \\
\Te &=& -\np\Pderive{\HABZ}{\te}(\phisABZ, \phirABZ, \te, \ke). \ase \label{eqn:ortho:te}
\ea
Notice that \eqnref{ortho:state-form}-\eqnref{ortho:algebraic} have the same structure as \eqref{eqn:appli:state-form}-\eqref{eqn:appli:algebraic}, i.e., corresponds to a Hamiltonian formulation, with~$(\usABZ - \Rs\isABZ, -\Rr\irABZ, -\Tl)$ as the vector of generalized forces. The $\abc$-to-$\ABZ$ transformation is what is called a canonical transformation, which preserves the Hamiltonian structure.

We emphasize that \eqnref{ortho:stator} and~\eqnref{ortho:rotor} automatically enforce the reciprocity conditions which must be satisfied by any reasonable motor model~\cite{MelkeW1990IAITo}; indeed,
\bano
	\frac{\partial\isA}{\partial\phisB} = \PPderive[\HABZ]{\phisB}{\phisA} &=& \PPderive[\HABZ]{\phisA}{\phisB} = \frac{\partial\isB}{\partial\phisA} \\
	\frac{\partial\irA}{\partial\phirB} = \PPderive[\HABZ]{\phirB}{\phirA} &=& \PPderive[\HABZ]{\phirA}{\phirB} = \frac{\partial\irB}{\partial\phirA}.
\eano

\subsection{$\ABZ$-to-$\dqZ$~transformation}\label{sec:dq}
Define the $\dqZ$~variables
\bano
\phisdqZ &:=& \Rot3(-\ts)\phisABZ \\
\phirdqZ &:=& \Rot3(\te-\ts)\phirABZ
\eano
and the scalar function
\bano
& & \HdqZ(\phisdqZ, \phirdqZ, \te, \ke ;\ts) \ane \\
& & \qquad := \HABZ\bigl(\Rot3(\ts)\phisdqZ, \Rot3(\ts - \te)\phirdqZ, \te, \ke\bigr).
\eano
We then find
\bano
\Tderive{\phisdqZ} &=& \Tderive{} \bigl(\Rot3(-\ts)\phisABZ\bigr) \\
&=& \Rot3(-\ts) \usABZ - \Rot3(-\ts)\Rs\Pderive{\HABZ}{\phisABZ} \\
& & \qquad -\>\ws\JJ3\Rot3(-\ts)\phisABZ \\
&=& \Rot3(-\ts) \usABZ - \Rs\Rot3(-\ts)\Rot3(\ts)\Pderive{\HdqZ}{\phisdqZ} \\
& & \qquad -\>\ws\JJ3\phisdqZ \\
\Tderive{\phirdqZ} &=& \Tderive{} \left(\Rot3(\te - \ts)\phirABZ\right) \\
&=& -\Rot3(\te - \ts) \Rr\Pderive{\HABZ}{\phirABZ} \\
& & \qquad -\>(\ws - \we)\JJ3\Rot3(\te - \ts)\phirABZ \\
&=& -\Rr\Rot3(\te - \ts)\Rot3(\ts - \te) \Pderive{\HdqZ}{\phirdqZ} \\
& & \qquad -\>(\ws - \we)\JJ3\phirdqZ; \\
\eano
moreover $\Te = \np\Transpose{\phirdqZ}\!\!\!\JJ3\irdqZ - \np\Pderive{\HdqZ}{\te}$, since
\bano
\Pderive{\HdqZ}{\te} &=& \Pderive{\HABZ}{\te} + \Jac{\HABZ}{\phirABZ} \Pderive{}{\te}\left(\Rot3(\ts - \te)\phirdqZ\right) \\
&=& \Pderive{\HABZ}{\te} - \Jac{\HABZ}{\phirABZ} {\Rot3}(\ts - \te) \JJ3 \phirdqZ \\
&=& \Pderive{\HABZ}{\te} + \Transpose{\phirdqZ}\!\!\!\JJ3\Pderive{\HdqZ}{\phirdqZ}.
\eano

Finally, the transformed system reads
\ba\label{eqn:dqZ:state-form}
\Tderive{\phisdqZ} &=& \usdqZ - \Rs\isdqZ - \JJ3\ws\phisdqZ \ase \\
\Tderive{\phirdqZ} &=& -\Rr\irdqZ - \JJ3(\ws - \we)\phirdqZ \ase \\
\Tderive{\te} &=& \we \ase \\
\frac{1}{\np} \Tderive{\ke} &=& \Te - \Tl, \ase
\ea
where $\usdqZ:=\Rot3(-\ts)\phisABZ$ and
\ba\label{eqn:dqZ:algebraic}
\isdqZ &:=& \Rot3(-\ts)\isABZ\ane\\ 
&=& \Pderive{\HdqZ}{\phisdqZ}(\phisdqZ, \phirdqZ, \te, \ke;\ts) \ase \label{eqn:dqZ:phii:stator} \\
\irdqZ &:=& \Rot3(\te-\ts)\irABZ\ane\\
&=& \Pderive{\HdqZ}{\phirdqZ}(\phisdqZ, \phirdqZ, \te, \ke;\ts) \ase \label{eqn:dqZ:phii:rotor} \ans \\
\we &=& \np^2\Pderive{\HdqZ}{\ke}(\phisdqZ, \phirdqZ, \te, \ke;\ts) \ase \label{eqn:dqZ:w} \\
\Te &=& -\np\Pderive{\HdqZ}{\te}(\phisdqZ, \phirdqZ, \te, \ke;\ts) + \np\Transpose{\phirdqZ}\!\!\!\JJ3\irdqZ. \ase\ans \label{eqn:dqZ:te}
\ea
The $\ABZ$-to-$\dqZ$ transformation is not canonical, as \eqnref{dqZ:state-form}-\eqnref{dqZ:algebraic} does not correspond to a Hamiltonian formulation, because of the new terms in the equations of the flux linkages and torque. This is not a problem for control purposes, and we can still consider~$\HdqZ$ as the ``pseudo-Hamiltonian'' of the system. Nevertheless, the reciprocity conditions
\begin{IEEEeqnarray*}{C'C}
	\frac{\partial\isd}{\partial\phisq}=\frac{\partial\isq}{\partial\phisd},\quad
	& \frac{\partial\ird}{\partial\phirq}=\frac{\partial\irq}{\partial\phird}
\end{IEEEeqnarray*}
are still satisfied.

\subsection{$\ABZ$-to-$\DQZ$~transformation}\label{sec:DQ}
Define the $\DQZ$~variables by
\bano
\phisDQZ &:=& \Rot3(-\te)\phisABZ \\
\phirDQZ &:=& \phirABZ
\eano
and the scalar function
\beno
\HDQZ(\phisDQZ, \phirDQZ, \te, \ke) := \HABZ(\Rot3(\te)\phisDQZ, \phirDQZ, \te, \ke).
\eeno
Proceeding as in the previous section, the transformed system reads
\ba\label{eqn:DQZ:state-form}
\Tderive{\phisDQZ} &=& \usDQZ - \Rs\isDQZ - \JJ3\we \phisDQZ \ase \\
\Tderive{\phirDQZ} &=& -\Rr\irDQZ \ase \\
\Tderive{\te} &=& \we \ase \\
\frac{1}{\np} \Tderive{\ke} &=& \Te - \Tl, \ase
\ea
where $\usDQZ:=\Rot3(-\te)\usABZ$ and
\ba\label{eqn:DQZ:algebraic}
\isDQZ &:=& \Rot3(-\te)\isABZ = \Pderive{\HDQZ}{\!\!\phisDQZ}(\phisDQZ\!\!, \phirDQZ\!\!, \te, \ke) \ans\ase \label{eqn:DQZ:current:stator} \\
\irDQZ &:=& \irABZ = \Pderive{\HDQZ}{\!\!\phirDQZ}(\phisDQZ, \phirDQZ, \te, \ke) \ase \label{eqn:DQZ:current:rotor} \\
\we &=& \np^2 \Pderive{\HDQZ}{\ke}(\phisDQZ, \phirDQZ, \te, \ke) \ase \label{eqn:DQZ:w}\\
\Te &=& -\np\Pderive{\HDQZ}{\te}(\phisDQZ, \phirDQZ, \te, \ke) \ane \\
& & -\>\np\Transpose{\phisDQZ}\!\!\!\JJ3\isDQZ. \ase \ans \label{eqn:DQZ:te}
\ea
The $\ABZ$-to-$\DQZ$ transformation is also not canonical, but the reciprocity conditions
\begin{IEEEeqnarray*}{C'C}
	\frac{\partial\isD}{\partial\phisQ}=\frac{\partial\isQ}{\partial\phisD},\quad
	& \frac{\partial\irD}{\partial\phirQ}=\frac{\partial\irQ}{\partial\phirD}
\end{IEEEeqnarray*}
are once again satisfied.

\subsection{Summary: equations in a general $xyz$ frame}
We can summarize the models
obtained in the frames~$\abc$, $\ABZ$, $\dqZ$ and~$\DQZ$ by writing them in a general $\xyz$~frame. The state equations~\eqnref{appli:state-form}, \eqnref{ortho:state-form}, \eqnref{dqZ:state-form} and~\eqnref{DQZ:state-form} then read
\ba\label{eqn:all:state-form}
\Tderive{\phisxyz} &=& \usxyz - \Rs\isxyz - \JJ3\Omega^{xyz}_s\phisxyz \ase \\
\Tderive{\phirxyz} &=& -\Rr\irxyz - \JJ3\Omega^{xyz}_r\phirxyz \ase \\
\Tderive{\te} &=& \we \ase \\
\frac{1}{\np} \Tderive{\ke} &=& \Te - \Tl, \ase
\ea
while the constitutive relations \eqnref{appli:algebraic}, \eqnref{ortho:algebraic}, \eqnref{dqZ:algebraic} and~\eqnref{DQZ:algebraic} read
\ba\label{eqn:all:algebraic}
\isxyz &=& \Pderive{\Hxyz}{\phisxyz}(\phisxyz, \phirxyz, \te, \ke) \ase \label{eqn:all:stator} \\
\irxyz &=& \Pderive{\Hxyz}{\phirxyz}(\phisxyz, \phirxyz, \te, \ke) \ase \label{eqn:all:rotor} \\
\we &=& \np^2\Pderive{\Hxyz}{\ke}(\phisxyz, \phirxyz, \te, \ke) \ase \\
\Te &=& -\np\Pderive{\Hxyz}{\te}(\phisxyz, \phirxyz, \te, \ke)+\Texyz. \ase
\ea
In~\eqnref{all:state-form}-\eqnref{all:algebraic}, the frame-specific terms are characterized by
\ba\label{eqn:all:fst}
(\Omega_s^{\abc}, \Omega_r^{\abc}, \Teabc) &=& (0, 0 , 0) \ase \\
(\Omega_s^{\ABZ}, \Omega_r^{\ABZ}, \TeABZ) &=& (0, 0 , 0) \ase \\
(\Omega_s^{\dqZ}, \Omega_r^{\dqZ}, \TedqZ) &=& (\ws, \ws - \we , \np\phirdqZ\JJ3\irdqZ) \ans \ase \\
(\Omega_s^{\DQZ}, \Omega_r^{\DQZ}, \TeDQZ) &=& (\we, 0 , -\np\phisDQZ\JJ3\isDQZ). \ase
\ea
These frame-specific terms are not zero in the $\dqZ$~and $\DQZ$~frames, because the transformations to these frames are not canonical, hence \eqnref{all:state-form}-\eqnref{all:algebraic}~is in these cases only a ``pseudo-Hamiltonian'' formulation.

%% file: symmetries.tex
\section{Using geometric symmetries of AC motors}\label{sec:symmetries}
The geometric symmetries enjoyed by any well-built motor constrain the form of its Hamiltonian, which facilitates its determination. The idea is simple: when a symmetry leaves the motor globally unchanged (meaning the new geometric configuration is indistinguishable from the original one), then the Hamiltonian is also unchanged.

\subsection{Stator symmetry: circular permutation of stator phases} \label{sec:symmetries:stator:perm}
If the three stator phases, which are identical, are circularly permuted, and the rotor is rotated by~$\frac{2\pi}{3}$, the motor is globally unchanged: indeed, rotating it by~$-\frac{2\pi}{3}$ gives back the initial configuration. The Hamiltonian must therefore satisfy
\beno
\Habc(\phisabc, \phirabc, \te, \ke) = \Habc\bigl(\Pabc\phisabc, \phirabc, \te + \tfrac{2\pi}{3}, \ke\bigr),
\eeno
where the phase permutation matrix~$\Pabc$ is defined as
\beno
\Pabc := \begin{pmatrix}
	0 & 0 & 1 \\
	1 & 0 & 0 \\
	0 & 1 & 0 
\end{pmatrix}.
\eeno
In the $\ABZ$~frame, this invariance property reads
\bano
& & \HABZ\bigl(\phisABZ, \phirABZ, \te, \ke\bigr) \\
& & \qquad = \Habc\bigl(\transpose{\CC}\phisABZ, \transpose{\CC}\phirABZ, \te, \ke\bigr) \\
& & \qquad = \Habc\bigl(\Pabc\transpose{\CC}\phisABZ, \transpose{\CC}\phirABZ, \te + \tfrac{2\pi}{3}, \ke\bigr) \\
& & \qquad = \HABZ\bigl(\CC\Pabc\transpose{\CC}\phisABZ, \CC\transpose{\CC}\phirABZ, \te + \tfrac{2\pi}{3}, \ke\bigr) \\
& & \qquad = \HABZ\bigl(\Rot3(\tfrac{2\pi}{3})\phisABZ, \phirABZ, \te + \tfrac{2\pi}{3}, \ke\bigr).
\eano
In the $\DQZ$~frame, this property is particularly simple,
\ba
& & \HDQZ\bigl(\phisDQZ, \phirDQZ, \te, \ke\bigr) \ane \\
& & \qquad = \HABZ\bigl(\Rot3(\te)\phisDQZ, \phirDQZ, \te, \ke\bigr) \ane \\
& & \qquad = \HABZ\bigl(\Rot3(\tfrac{2\pi}{3})\Rot3(\te)\phisDQZ, \phirDQZ, \te + \tfrac{2\pi}{3}, \ke\bigr) \ane \\
& & \qquad = \HDQZ\bigl(\Rot3(-\te - \tfrac{2\pi}{3})\Rot3(\tfrac{2\pi}{3})\Rot3(\te)\phisDQZ, \ane \\
& & \qquad\qquad \phirDQZ, \te + \tfrac{2\pi}{3}, \ke\bigr) \ane \\
& & \qquad = \HDQZ\bigl(\phisDQZ, \phirDQZ, \te + \tfrac{2\pi}{3}, \ke\bigr); \label{eqn:symmetries:stator:perm:DQZ}
\ea
in other words, the Hamiltonian in the $\DQZ$~frame is $\frac{2\pi}{3}$\nobreakdash-periodic with respect to $\te$.
Finally, in the $\dqZ$~frame,
\ba
& & \HdqZ\bigl(\phisdqZ, \phirdqZ, \te, \ke\bigr) \ane \\
& & \qquad = \HABZ\bigl(\Rot3(\ts)\phisdqZ, \Rot3(\ts - \te)\phirdqZ, \te, \ke\bigr) \ane \\
& & \qquad = \HABZ\bigl(\Rot3(\tfrac{2\pi}{3})\Rot3(\ts)\phisdqZ, \ane \\
& & \qquad\qquad \Rot3(\ts - \te)\phirdqZ, \te + \tfrac{2\pi}{3}, \ke\bigr) \ane \\
& & \qquad = \HdqZ\bigl(\Rot3(-\ts)\Rot3(\tfrac{2\pi}{3})\Rot3(\ts)\phisdqZ, \ane \\
& & \qquad\qquad \Rot3(\te + \tfrac{2\pi}{3} - \ts)\Rot3(\ts - \te)\phirdqZ, \te + \tfrac{2\pi}{3}, \ke\bigr) \ane \\
& & \qquad = \HdqZ\bigl(\Rot3(\tfrac{2\pi}{3})\phisdqZ, \Rot3(\tfrac{2\pi}{3})\phirdqZ, \te + \tfrac{2\pi}{3}, \ke\bigr). \label{eqn:symmetries:stator:perm:dqZ}
\ea
	
\subsection{Stator symmetry: reversal of stator flux linkages} \label{sec:symmetries:stator:rev}
As the opposite poles of the stator are identical, if the stator flux linkages are reversed and the rotor is rotated by~$\pi$, the motor is globally unchanged: indeed, rotating it by~$-\pi$ gives back the initial configuration. The Hamiltonian must therefore satisfy
\beno
\Habc(\phisabc, \phirabc, \te, \ke) = \Habc(-\phisabc, \phirabc, \te + \pi, \ke).
\eeno
Proceeding as in the previous section, this yields
\ba
& & \HABZ(\phisABZ, \phirABZ, \te, \ke) \ane \\
& & \qquad = \HABZ(-\phisABZ, \phirABZ, \te + \pi, \ke) \ane \\
& & \HDQZ(\phisD, \phisQ, \phisZ, \phirD, \phirQ, \phirZ, \te, \ke) \ane \\
& & \qquad = \HDQZ(\phisD, \phisQ, -\phisZ, \phirD, \phirQ, \phirZ, \te + \pi, \ke) \label{eqn:symmetries:stator:rev:DQZ} \\
& & \HdqZ(\phisd, \phisq, \phisZ, \phird, \phirq, \phirZ, \te, \ke) \ane \\
& & \qquad = \HdqZ(-\phisd, -\phisq, -\phisZ, -\phird, -\phirq, \phirZ, \te + \pi, \ke). \ans\label{eqn:symmetries:stator:rev:dqZ}
\ea

\subsection{Rotor symmetry: reversal of rotor flux linkages} \label{sec:symmetries:rotor:rev} 
As the opposite poles of the rotor are identical, if the rotor flux linkages are reversed and the rotor is rotated by~$\pi$, the motor is globally unchanged: indeed, rotating it by~$-\pi$ gives back the initial configuration. The Hamiltonian must therefore satisfy
\beno
\Habc(\phisabc, \phirabc, \te, \ke) = \Habc(\phisabc, -\phirabc, \te + \pi, \ke).
\eeno
Proceeding as in the previous section, this yields
\ba
& & \HABZ(\phisABZ, \phirABZ, \te, \ke) \ane \\
& & \qquad = \HABZ(\phisABZ, -\phirABZ, \te + \pi, \ke) \ane \\
& & \HDQZ(\phisD, \phisQ, \phisZ, \phirD, \phirQ, \phirZ, \te, \ke) \ane \\
& & \qquad = \HDQZ(-\phisD\!, -\phisQ\!, \phisZ, -\phirD\!, -\phirQ\!, -\phirZ, \te + \pi, \ke)  \ans\label{eqn:symmetries:rotor:rev:DQZ} \\
& & \HdqZ(\phisd, \phisq, \phisZ, \phirq, \phirq, \phirZ, \te, \ke) \ane \\
& & \qquad = \HdqZ(\phisd, \phisq, \phisZ, \phird, \phirq, -\phirZ, \te + \pi, \ke). \label{eqn:symmetries:rotor:rev:dqZ}
\ea

\subsection{Rotor symmetry: circular permutation of rotor phases/bars} \label{sec:symmetries:rotor:perm}
In an induction motor, the rotor is invariant by a rotation of angle~$\eta$; for a wound-rotor motor $\eta = \frac{2\pi}{3}$, since the three phases are identical; for a squirrel-cage motor $\eta = \frac{2\pi}{\nr}$ (with $\nr$~the number of bars per pole pair), since the bars are identical. Though this is a symmetry in the $\abc$~frame, it is easier to express it in the $\ABZ$~frame: if the rotor flux linkage is rotated by~$-\eta$ around the $0$\nobreakdash-axis and the rotor is rotated by~$\eta$, the motor is globally unchanged. The Hamiltonian must therefore satisfy
\beno
\HABZ(\phisABZ\!, \phirABZ\!, \te, \ke) = \HABZ\bigl(\phisABZ\!, \Rot3(-\eta)\phirABZ\!, \te + \eta, \ke\bigr).
\eeno
In the $\dqZ$~frame,
\ba
& & \HdqZ\bigl(\phisdqZ, \phirdqZ, \te, \ke\bigr) \ane \\
& & \qquad = \HABZ\bigl(\Rot3(\ts)\phisdqZ, \Rot3(\ts - \te)\phirdqZ, \te + \eta, \ke\bigr) \ane \\
& & \qquad = \HABZ\bigl(\Rot3(\ts)\phisdqZ, \Rot3(-\eta)\Rot3(\ts - \te)\phirdqZ, \te + \eta, \ke\bigr) \ane \\
& & \qquad = \HdqZ\bigl(\Rot3(-\ts)\Rot3(\ts)\phisdqZ, \ane \\
& & \qquad\qquad\Rot3(\te + \eta - \ts)\Rot3(-\eta)\Rot3(\ts - \te)\phirdqZ, \te + \eta, \ke\bigr) \ane \\
& & \qquad = \HdqZ(\phisdqZ, \phirdqZ, \te + \eta, \ke). \label{eqn:symmetries:rotor:perm:dqZ}
\ea

\subsection{Stator and rotor symmetry: exchange of phases} \label{sec:symmetries:swap}
When the rotor is symmetric with respect to a plane (without loss of generality the $a$\nobreakdash-axis belongs to this plane), if the $b$~and $c$~phases are swapped in both the stator and the rotor, and the direction of the rotation is reversed, the motor is globally unchanged: indeed, a reflection gives back the initial configuration. The Hamiltonian must therefore satisfy
\beno
\Habc(\phisabc, \phirabc, \te, \ke) = \Habc(\Oabc\phisabc, \Oabc\phirabc, -\te, -\ke),
\eeno
where
\beno
\Oabc := \begin{pmatrix}
	1 & 0 & 0 \\
	0 & 0 & 1 \\
	0 & 1 & 0 
\end{pmatrix}.
\eeno
Expressed in the $\ABZ$~frame, this reads
\bano
& & \HABZ(\phisABZ, \phirABZ, \te, \ke) \ane \\
& & \qquad = \HABZ(\OABZ\phisABZ, \OABZ\phirABZ, -\te, -\ke),
\eano
where
\beno
\OABZ := \CC\Oabc\CC^{-1} = \begin{pmatrix}
	1 & 0 & 0 \\
	0 & -1 & 0 \\
	0 & 0 & 1
\end{pmatrix}.
\eeno
In the $\DQZ$ and $\dqZ$~frames,
\bano
& & \HDQZ\bigl(\phisDQZ, \phirDQZ, \te, \ke\bigr) \\
& & \qquad = \HABZ\bigl(\Rot3(\te)\phisDQZ, \phirDQZ, \te, \ke\bigr) \\
& & \qquad = \HABZ\bigl(\OABZ\Rot3(\te)\phisDQZ, \OABZ\phirDQZ, -\te, -\ke\bigr) \\
& & \qquad = \HDQZ\bigl(\Rot3(\te)\OABZ\Rot3(\te)\phisDQZ, \\
& & \qquad\qquad \OABZ\phirDQZ, -\te, -\ke\bigr) \\
& & \HdqZ\bigl(\phisdqZ, \phirdqZ, \te, \ke\bigr) \\
& & \qquad = \HdqZ\bigl(\Rot3(\ts)\OABZ\Rot3(\ts)\phisdqZ, \\
& & \qquad\qquad \Rot3(\ts - \te)\OABZ\Rot3(\ts - \te)\phirdqZ, -\te, -\ke\bigr). \\
\eano
Rewritten component-wise, this yields the parity conditions
\ba
& & \HDQZ(\phisD, \phisQ, \phisZ, \phirD, \phirQ, \phirZ, \te, \ke) \ane\\
& & \qquad = \HDQZ(\phisD, -\phisQ, \phisZ, \phirD, -\phirQ, \phirZ, -\te, -\ke) \label{eqn:symmetries:rotor:swapQ:DQZ}  \\
& & \HdqZ(\phisd, \phisq, \phisZ, \phird, \phirq, \phirZ, \te, \ke) \ane  \\
& & \qquad = \HdqZ(\phisd, -\phisq, \phisZ, \phird, -\phirq, \phirZ, -\te, -\ke).\label{eqn:symmetries:rotor:swapQ:dqZ}
\ea

If the rotor is moreover symmetric with respect to the orthogonal plane, the same reasoning yields
\ba
& & \HDQZ(\phisD, \phisQ, \phisZ, \phirD, \phirQ, \phirZ, \te, \ke) \ane \\
& & \qquad = \HDQZ(-\phisD, \phisQ, \phisZ, -\phirD, \phirQ, \phirZ, -\te, -\ke) \label{eqn:symmetries:rotor:swapD:DQZ} \\
& & \HdqZ(\phisd, \phisq, \phisZ, \phird, \phirq, \phirZ, \te, \ke) \ane \\
& & \qquad = \HdqZ(-\phisd, \phisq, \phisZ, -\phird, \phirq, \phirZ, -\te, -\ke). \label{eqn:symmetries:rotor:swapD:dqZ}
\ea

%% file: connections.tex
\section{Using connection constraints}\label{sec:connections}
Up to now, we have been considering AC motors with unconnected windings. In this section, we take these connections into account, which imposes constraints on the currents involved in the equations derived in \secref{appli}. We first study how to handle constraints in the general framework of analytical mechanics and then apply this to AC motors.

\subsection{Analytical mechanics with constraints (see e.g.~\cite{Lurie2002book})}\label{sec:connection:constraints}
Constraints in analytical mechanics are of two kinds: holonomic constraints $g(q)=0$, where $g(q)\in\mathbb{R}^m$; nonholonomic constraints $K(q)\tderive{q}-k(q)=0$, where $K(q)$~is a full-rank $m\times n$~matrix and~$k(q)\in\mathbb{R}^m$. Notice that in view of~\eqnref{meca:dynamic:qdot} the holonomic constraint $g(q)=0$ naturally gives rise to the nonholonomic constraint $g'(q)\tderive{q}=0$; this  nonholonomic constraint is equivalent to the original holonomic constraint, up to a constant which is determined by the initial value $q(0)$. The holonomic constraint $g(q)=0$ can thus be handled in two (equivalent) ways: on the one hand, assuming it is explicitly solvable as $q_2=h(q_1)$, we can define the reduced Lagrangian
\beno
\Lg_r(q_1,\tderive{q_1}) := \Lg\bigl(q_1,h(q_1),\tderive{q_1},h'(q_1)\tderive{q_1}\bigr),
\eeno
and proceed exactly as in section~\ref{sec:basics} with $\Lg_r$ instead of~$\Lg$; on the other hand, rather than explicitly solving the constraint, which may be cumbersome in practice, we can consider only the associated nonholonomic constraint $g'(q)\tderive{q}=0$, and handle it as any other nonholonomic constraint.

The nonholonomic constraint $K(q)\tderive{q}-k(q)=0$ is handled by introducing the so-called Lagrange multiplier~$\mu \in \mathbb{R}^m$ and using the constrained version of the Euler-Lagrange equations
\ba\label{eqn:connection:euler-lagrange}
\Tderive{q} &=& \tderive{q} \ase \label{eqn:connection:euler-lagrange:qdot} \\
\Tderive{} \Pderive{\Lg}{\tderive{q}} &=& \Pderive{\Lg}{q}(q, \tderive{q}) + \mathcal{F} + \transpose{K}(q)\mu \ase \label{eqn:connection:euler-lagrange:q} \\
0 &=& K(q)\tderive{q} - k(q). \ase \label{eqn:connection:euler-lagrange:lambda}
\ea
To obtain the Hamiltonian formulation, we then proceed as in \secref{basics}, but use~\eqnref{connection:euler-lagrange:qdot}-\eqnref{connection:euler-lagrange:q} instead of~\eqnref{meca:dynamic}. This yields
\ba\label{eqn:connection:da-system}
\Tderive{p} &=& -\Pderive{\Hm}{q}(p, q) + \mathcal{F} + \transpose{K}(q)\mu \ase \label{eqn:connection:da-system:p} \\
\Tderive{q} &=& \Pderive{\Hm}{p}(p, q) \ase \label{eqn:connection:da-system:q} \\
0 &=& K(q)\Pderive{\Hm}{p}(p, q) - k(q), \ase \label{eqn:connection:da-system:lambda}
\ea
where we have replaced~$\tderive{q}$ in~\eqnref{connection:euler-lagrange:lambda} by its expression~\eqnref{connection:da-system:q}. This is a differential-algebraic system with $2n$~differential equations and $m$~algebraic constraints for the $2n + m$~unknowns $p$, $q$ and~$\mu$. Of course, the initial condition~$(p_0,q_0)$ must respect the constraint~\eqnref{connection:da-system:lambda}. 

Such a differential-algebraic may be difficult to handle, e.g. for simulation purposes. If needed, \eqnref{connection:da-system} may be rewritten as a purely differential system: first differentiate the constraint~\eqnref{connection:da-system:lambda} and use~\eqnref{connection:da-system:p}-\eqnref{connection:da-system:q}, which yields
\bano
0 &=&  K \Bigl(\pderive{}{p} \Pderive{\Hm}{p} \Tderive{p} + \pderive{}{q} \Pderive{\Hm}{p} \Tderive{q}\Bigr) + \Bigl(K' \Tderive{q}\Bigr) \Pderive{\Hm}{p} - k'\Tderive{q} \\
&=& K \pderive{}{p} \Pderive{\Hm}{p}\transpose{K}\mu + K \pderive{}{p} \Pderive{\Hm}{p}(\mathcal{F} - \Pderive{\Hm}{q}) \\
& & \qquad + \bigr(K \pderive{}{q} \Pderive{\Hm}{p} + K' \Pderive{\Hm}{p} - k'\bigl)\Pderive{\Hm}{p};
\eano
since $K$ is full rank and~$\pderive{}{p} \Pderive{\Hm}{p}$ is invertible by assumption, $K \pderive{}{p}\Pderive{\Hm}{p} \transpose{K}$~is invertible and the previous equation can be solved for the Lagrange multiplier as
\be\label{eqn:connection:lambda:H}
\mu = \mathcal{M}(p, q, \mathcal{F}).
\ee
The purely differential version of~\eqnref{connection:da-system} is then
\ba\label{eqn:connection:state-form}
\Tderive{p} &=& -\Pderive{\Hm}{q}(p, q) + \mathcal{F} + \transpose{K}(q)\mathcal{M}(p, q, \mathcal{F}) \ase \label{eqn:connection:state-form:p} \\
\Tderive{q} &=& \Pderive{\Hm}{p}(p, q), \ase \label{eqn:connection:state-form:q} 
\ea
which can be readily simulated starting from any initial condition~$(p_0, q_0)$ such that~$K(q_0)\Pderive{\Hm}{p}(p_0, q_0) - k(q_0) = 0$.			
			
\subsection{Short-circuited rotor windings}\label{sec:connection:rotor}
\begin{figure}
	\centering
	\subfloat[Short-circuited rotor windings.\label{fig:connection:short-rotor}]{\includegraphics[width=120pt]{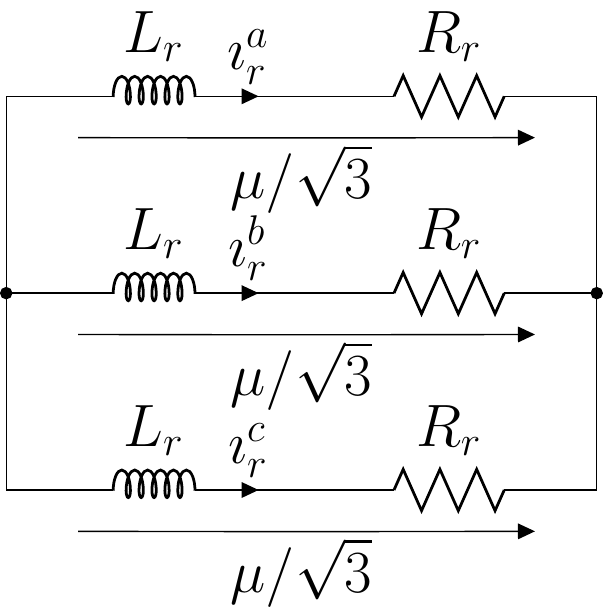}} \hfill
	\subfloat[Fictitious windings.\label{fig:connection:norotor}]{\includegraphics[width=120pt]{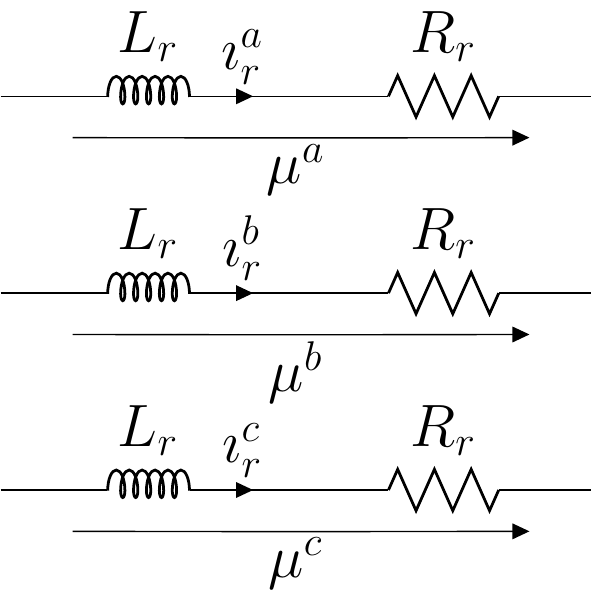}}
	\caption{Main rotor connection schemes.}
	\label{fig:connection:rotor}
\end{figure}

The rotor windings of AC motors are usually short-circuited, see \figref{connection:short-rotor}, which imposes the constraint $\ira + \irb + \irc = 0$.
In any frame~$\xyZ$, this constraint reads
\be\label{eqn:connection:rotor:constraint}
0 = \irZ := \Pderive{\HxyZ}{\phirZ}(\phisxyZ, \phirxy, \phirZ, \te, \ke);
\ee
if the (pseudo)-Hamiltonian is not degenerate, then $\Jac{}{\phirZ}\Pderive{\HxyZ}{\phirZ} \neq 0$, and $\phirZ = \Lambda_r^0(\phisxyZ, \phirxy, \te, \ke)$  by the implicit function theorem. Define now the function
\bano
& & \HcxyZ(\phisxyZ, \phirxy, \te, \ke) \ane \\
& & \qquad:= \HxyZ\bigl(\phisxyZ, \phirxy, \Lambda_r^0(\phisxyZ, \phirxy, \te, \ke), \te, \ke\bigr).
\eano
Its gradients with respect to $\chi \in \{\phisxyZ, \phirxy, \te, \ke\}$ are
\bano
& & \Pderive{\HcxyZ}{\chi}(\phisxyZ, \phirxy, \te, \ke) \\
& & \quad = \Pderive{\HxyZ}{\chi}(\ldots) + \Pderive{\Lambda_r^0}{\chi}(\ldots)\Pderive{\HxyZ}{\phirZ}(\ldots) \\
& & \quad = \Pderive{\HxyZ}{\chi}\bigl(\phisxyZ, \phirxy, \Lambda_r^0(\phisxyZ, \phirxy, \te, \ke), \te, \ke\bigr),
\eano
where $\Pderive{\HxyZ}{\phirZ}(\ldots) = 0$ by~\eqnref{connection:rotor:constraint}. The (pseudo)-Hamiltonian equations~\eqnref{all:state-form}-\eqnref{all:algebraic} with constraint~\eqnref{connection:rotor:constraint} then split into two decoupled subsystems:
\begin{itemize}
	\item the unconstrained $\xy$-subsystem, with state equations
	\ba
	\Tderive{\phisxyZ} &=& \usxyZ - \Rs\isxyZ - \JJ3\Omega^{xy0}_s\phisxyZ \ase \\
	\Tderive{\phirxy} &=& -\Rr\irxy - \JJ2\Omega^{xy0}_r\phirxy \ase \\
	\Tderive{\te} &=& \we \ase \\
	\frac{1}{\np} \Tderive{\ke} &=& \Te - \Tl, \ase
	\ea
	and constitutive relations
	\ba\label{eqn:xyZ:state-formWound}
	\isxyZ &=& \Pderive{\HcxyZ}{\phisxyZ}(\phisxyZ, \phirxy, \te, \ke) \ase \\
	\irxy &=& \Pderive{\HcxyZ}{\phirxy}(\phisxyZ, \phirxy, \te, \ke) \ase \\
	\we &=& \np^2\Pderive{\HcxyZ}{\ke}(\phisxyZ, \phirxy, \te, \ke) \ase \\
	\Te &=& -\np\Pderive{\HcxyZ}{\te}(\phisxyZ, \phirxy, \te, \ke) +\TexyZ; \ase
	\ea
	notice $\TedqZ=\np\phirdqZ\JJ3\irdqZ=\np\phirdq\JJ2\irdq$ by~\eqnref{all:state-form}
	\item the constrained $0$-subsystem
	\ba\label{eqn:xyZ:algebraicWound}
	\Tderive{\phirZ} &=& \mu \ase \\
	0 &=& \Pderive{\HxyZ}{\phirZ}(\phisxyZ, \phirxy, \phirZ, \te, \ke). \ase
	\ea
	For simulation purposes, $\mu$~can be computed as a function of the state variables using~\eqnref{connection:lambda:H}.
	Notice $\mu/\sqrt{3}$~can be seen as the voltage drop across the short-circuited windings~$a$, $b$ and~$c$, see \figref{connection:short-rotor}. 
\end{itemize}

\subsection{No rotor windings}\label{sec:connection:norotor}
PMSMs and SynRMs do not have rotor windings, which can be translated into the 3-dimensional constraint 
\beno
0 = \irxyz(\phisxyz, \phirxyz, \te, \ke) := \Pderive{\Hxyz}{\phirxyz}(\phisxyz, \phirxyz, \te, \ke)
\eeno
in any frame~$\xyz$; this simply means there are no currents in the (fictitious) rotor windings. Proceeding then exactly as in \secref{connection:rotor}, we can invert the constraint and define the function
\beno
\Hzxyz(\phisxyz, \te, \ke) := \Hxyz\bigl(\phisxyz, \Lambda_r^{xyz}(\phisxyz, \te, \ke), \te, \ke\bigr),
\eeno
yielding the two decoupled subsystems:
\begin{itemize}
	\item the unconstrained stator subsystem, with state equations
	\ba\label{eqn:xyZ:state-formNo}
	\Tderive{\phisxyz} &=& \usxyz - \Rs\isxyz - \JJ3\Omega^{xyz}_s\phisxyz \ase \\
	\Tderive{\te} &=& \we \ase \\
	\frac{1}{\np} \Tderive{\ke} &=& \Te - \Tl, \ase
	\ea
	and constitutive relations
	\ba\label{eqn:xyZ:algebraicNo}
	\isxyz &=& \Pderive{\Hzxyz}{\phisxyz}(\phisxyz, \te, \ke) \ase \\
	\we &=& \np^2\Pderive{\Hzxyz}{\ke}(\phisxyz, \te, \ke) \ase \\
	\Te &=& -\np\Pderive{\Hzxyz}{\te}(\phisxyz, \te, \ke) + \Texyz; \ase
	\ea
	notice $\TedqZ=\np\phirdqZ\JJ3\irdqZ=0$, see~\eqnref{all:state-form}
	\item the (3-dimensional) constrained rotor subsystem
	\ba
	\Tderive{\phirxyz} &=& \mu - \JJ3\Omega^{xyz}_r\phirxyz \ase \\
	0 &=& \Pderive{\Hxyz}{\phirxyz}(\phisxyz, \phirxyz, \te, \ke). \ase
	\ea
	For simulation purposes, $\mu$~can be computed as a function of the state variables using~\eqnref{connection:lambda:H}.
	Notice $\mu$~can be seen as the vector $(\mu^a, \mu^b, \mu^c)$ of voltage drops across the (fictitious) rotor windings~$a$, $b$ and~$c$, see \figref{connection:norotor}.
\end{itemize}

\subsection{Star connection scheme}\label{sec:connection:star}
\begin{figure}
	\centering
	\subfloat[Unconnected stator windings.\label{fig:connection:unconnected}]{\includegraphics[width=120pt]{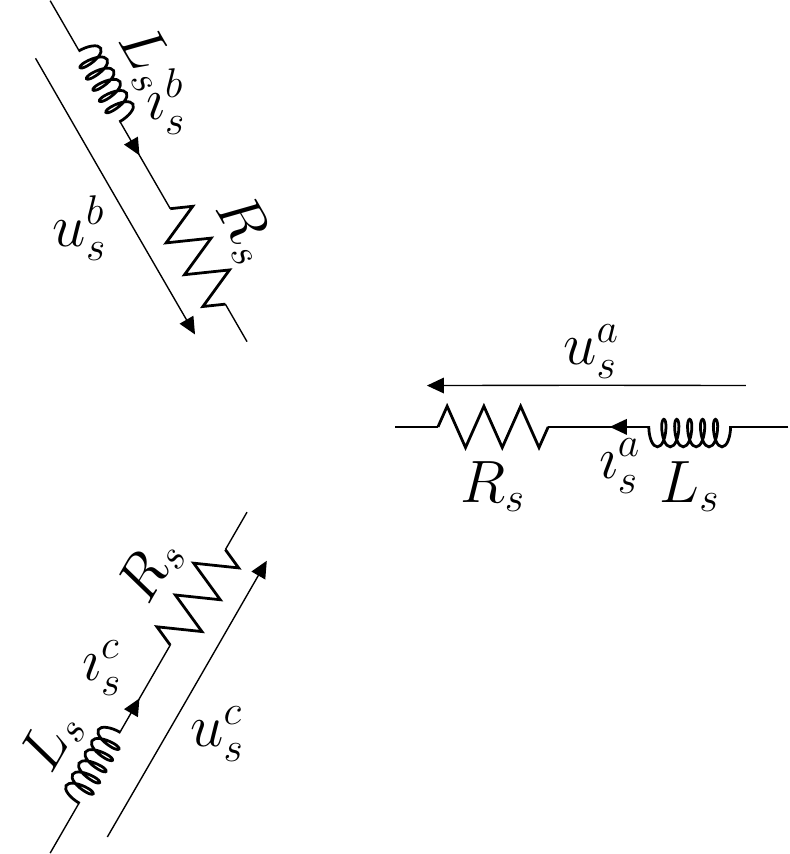}} \hfill
	\subfloat[The star connection scheme.\label{fig:connection:star}]{\includegraphics[width=120pt]{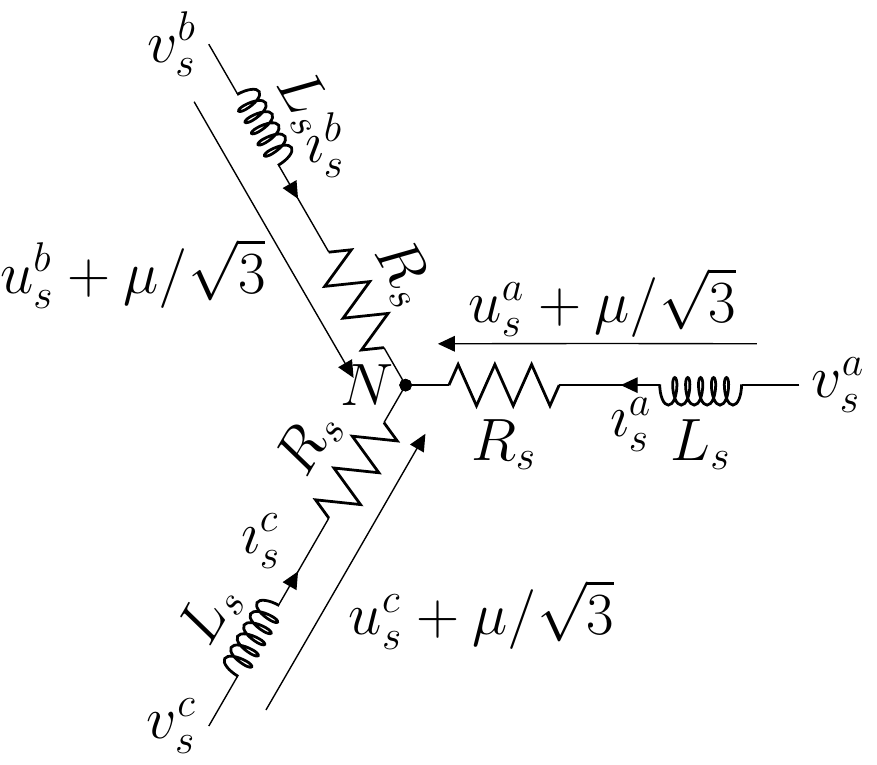}} \hfill
	\caption{Main stator connection schemes.}
	\label{fig:connection:stator}
\end{figure}

The stator windings of AC motors are usually star-connected, see~\figref{connection:star}, which imposes the constraint $\isa + \isb + \isc = 0$.
In this connection scheme, the potentials~$\vsa$, $\vsb$ and~$\vsc$ are impressed, and not the voltages drops across the windings.
The constrained version of \eqnref{appli:state-form:phis} reads
\beno
\Tderive{\phisabc} = \usabc - \Rs\isabc + \frac{1}{\sqrt{3}}\begin{pmatrix}
	1 \\
	1 \\
	1
\end{pmatrix}\mu.
\eeno
Notice the voltage drops across the windings are by definition
\bano
\vsa - v_\starpoint &:=& \Tderive{\phisa} + \Rs\isa = \usa + \mu/\sqrt{3} \\
\vsb - v_\starpoint &:=& \Tderive{\phisb} + \Rs\isb = \usb + \mu/\sqrt{3} \\
\vsc - v_\starpoint &:=& \Tderive{\phisc} + \Rs\isc = \usc + \mu/\sqrt{3}.
\eano
In any frame~$\xyZ$, this gives
\ba
\vsxy &=& \usxy \ase \\
\vsZ - \sqrt{3}v_\starpoint &=& \usZ + \mu, \ase \label{eqn:connection:star:vn}
\ea
and the constraint $\isa + \isb + \isc = 0$ reads
\beno
0 = \isZ := \Pderive{\HxyZ}{\phisZ}(\phisxy, \phisZ, \phirxyZ, \te, \ke).
\eeno
Proceeding then exactly as in \secref{connection:rotor}, we can invert the constraint and define the function
\ba
& & \HsxyZ(\phisxy, \phirxyZ, \te, \ke) \ane \\
& & \qquad := \HxyZ\bigl(\phisxy, \Lambda_s^0(\phisxy, \phirxyZ, \te, \ke), \phirxyZ, \te, \ke\bigr),
\ea
yielding the two decoupled subsystems:
\begin{itemize}
	\item the unconstrained $\xy$-subsystem, with state equations
	\ba\label{eqn:connection:star:state-form}
	\Tderive{\phisxy} &=& \vsxy - \Rs\isxy - \JJ2\Omega^{\xyZ}_s\phisxy \ase \\
	\Tderive{\phirxyZ} &=& -\Rr\irxyZ - \JJ3\Omega^{\xyZ}_r\phirxyZ \ase \\
	\Tderive{\te} &=& \we \ase \\
	\frac{1}{\np} \Tderive{\ke} &=& \Te - \Tl, \ase
	\ea
	and constitutive relations
	\ba\label{eqn:connection:star:algebraic}
	\isxy &=& \Pderive{\HsxyZ}{\phisxy}(\phisxy, \phirxyZ, \te, \ke) \ase \\
	\irxyZ &=& \Pderive{\HsxyZ}{\phirxyZ}(\phisxy, \phirxyZ, \te, \ke) \ase \\
	\we &=& \np^2\Pderive{\HsxyZ}{\ke}(\phisxy, \phirxyZ, \te, \ke) \ase \\
	\Te &=& -\np\Pderive{\HsxyZ}{\te}(\phisxy, \phirxyZ, \te, \ke) + \TexyZ; \ase
	\ea
	notice $\TeDQZ=-\np\phisDQ\JJ2\isDQ$, see~\eqnref{all:state-form}
	\item the constrained $0$-subsystem
	\ba\label{eqn:connection:star:zero}
	\Tderive{\phisZ} &=& \usZ + \mu \ase \label{eqn:connection:star:phiz} \\
	0 &=& \Pderive{\HxyZ}{\phisZ}(\phisxyZ, \phirxyZ, \te, \ke). \ase \label{eqn:connection:star:constraint}
	\ea
	Time differentiating \eqnref{connection:star:constraint} yields
	\beno
	0 = \Jac{}{\phisZ} \Pderive{\HxyZ}{\phisZ} \Tderive{\phisZ} + f(\phisxy, \phisZ, \phirxyZ, \te; \vsxy, \Tl),
	\eeno
	hence the purely differential version of~\eqnref{connection:star:zero} is
	\beno
	\Tderive{\phisZ} = -(\Jac{}{\phisZ} \Pderive{\HxyZ}{\phisZ} )^{-1}f(\phisxy, \phisZ, \phirxyZ, \te; \vsxy, \Tl);
	\eeno
	this is of course the same as using~\eqnref{connection:lambda:H}, but more direct. Since this equation is expressed in terms of the state variables and the (known) inputs~$\vsxy$ and~$\Tl$, it can be simulated. If desired, the potential~$v_\starpoint$ of the star point can then be computed using~\eqnref{connection:star:vn}.
\end{itemize}

\subsection{Summary: constraints, change of frames, and decoupling}\label{sec:connections:summary}
In the literature, it is usually taken for granted that control-oriented models of AC motors can be expressed in the $\DQ$ (or $\dq$) frame, i.e., involve a transformation from the physical $\abc$ frame and decoupling from the $0$-axis. Whereas this is easily justified for unsaturated sinusoidal models, the situation is more delicate for saturated models. Indeed, saturations physically takes place in the $\abc$ frame, and it is not obvious that the corresponding nonlinearities behave well under a change of frames and result in decoupling of the $0$-axis. One interest of the energy-based modeling approach is to fully justify this decoupling for star-connected motors.

For a star-connected motor with rotor windings, we can combine the results of sections \ref{sec:connection:rotor} and~\ref{sec:connection:star} and define a (pseudo)-Hamiltonian~$\Hcsxy(\phisxy, \phirxy, \te, \ke)$ independent of the $0$-variables. This yields the unconstrained $\xy$-subsystem, with state equations 
	\ba
	\Tderive{\phisxy} &=& \vsxy - \Rs\isxy - \JJ2\Omega^{\xy}_s\phisxy \ase \\
	\Tderive{\phirxy} &=& -\Rr\irxy - \JJ2\Omega^{xy}_r\phirxy \ase \\
	\Tderive{\te} &=& \we \ase \\
	\frac{1}{\np} \Tderive{\ke} &=& \Te - \Tl, \ase
	\ea
	and constitutive relations
	\ba
	\isxy &=& \Pderive{\Hcsxy}{\phisxy}(\phisxy, \phirxy, \te, \ke) \ase \\
	\irxy &=& \Pderive{\Hcsxy}{\phirxy}(\phisxy, \phirxy, \te, \ke) \ase \\
	\we &=& \np^2\Pderive{\Hcsxy}{\ke}(\phisxy, \phirxy, \te, \ke) \ase \\
	\Te &=& -\np\Pderive{\Hcsxy}{\te}(\phisxy, \phirxy, \te, \ke) +\Texy, \ase
	\ea
	where the frame-specific terms are given by
	\ba\label{eqn:all:fstxy}
	(\Omega_s^{\AB}, \Omega_r^{\AB}, \TeAB) &=& (0, 0 , 0) \ase \\
	(\Omega_s^{\dq}, \Omega_r^{\dq}, \Tedq) &=& (\ws, \ws - \we , \np\phirdq\JJ2\irdq) \ans \ase \\
	(\Omega_s^{\DQ}, \Omega_r^{\DQ}, \TeDQ) &=& (\we, 0 , -\np\phisDQ\JJ2\isDQ). \ase
	\ea
	
For a star-connected motor without rotor windings, we can combine the results of sections \ref{sec:connection:norotor} and~\ref{sec:connection:star} and define a (pseudo)-Hamiltonian~$\Hzsxy(\phisxy, \te, \ke)$ independent of the $0$-variables. This yields the unconstrained $\xy$-subsystem, with state equations 
\ba
\Tderive{\phisxy} &=& \vsxy - \Rs\isxy - \JJ2\Omega^{\xy}_s\phisxy \ase \\
\Tderive{\te} &=& \we \ase \\
\frac{1}{\np} \Tderive{\ke} &=& \Te - \Tl, \ase
\ea
and constitutive relations
\ba
\isxy &=& \Pderive{\Hzsxy}{\phisxy}(\phisxy, \te, \ke) \ase \\
\we &=& \np^2\Pderive{\Hzsxy}{\ke}(\phisxy, \te, \ke) \ase \\
\Te &=& -\np\Pderive{\Hzsxy}{\te}(\phisxy, \te, \ke) +\Texy, \ase
\ea
where the frame-specific terms are given by
\ba\label{eqn:all:fstxy:norotor}
(\Omega_s^{\AB}, \TeAB) &=& (0, 0) \ase \\
(\Omega_s^{\dq}, \Tedq) &=& (\ws, 0) \ase \\
(\Omega_s^{\DQ}, \TeDQ) &=& (\we, -\np\phisDQ\JJ2\isDQ). \ase
\ea

In both previous cases, the full (pseudo)-Hamiltonian $\HxyZ$ is required only if the stator and rotor $0$-subsystems are needed, which is usually not the case for control purposes.

%% file: simpleModels.tex
\section{Unsaturated sinusoidal models of AC motors}\label{sec:models}
In this section, we recover the classical models of the literature, namely unsaturated sinusoidal models. This is achieved by encoding in the (pseudo)-Hamiltonians the corresponding macroscopic assumptions, and taking advantage of the symmetries of \secref{symmetries} and the constraints of \secref{connections}. It is remarkable that no information beyond that is needed. The simplicity of the derivation is to be contrasted with the usual approach, see e.g.~\cite{Vas1998book,KrauseWSP2013book,Chiasson2005book}.

In electro-mechanical devices, the energy is of two kinds: mechanical energy and magnetic energy. As we are not interested in modeling the load, we consider here the simplest form for the mechanical energy, namely~$\frac{\ke^2}{2\Jl\np^2}$, which corresponds to a simple load with constant inertia; notice \eqnref{DQZ:w} reads in this case $\we=\frac{\ke}{J}$. On the other hand, experimental observations show that the magnetic energy does not depend on the rotor velocity, hence on the kinetic momentum. We therefore seek the (pseudo)-Hamiltonian in the form
\beno
\HDQZ(\phisDQZ\!\!, \phirDQZ, \te, \ke) = \frac{\ke^2}{2J\np^2} + \HmDQZ(\phisDQZ\!\!, \phirDQZ\!\!, \te),
\eeno
where $\HmDQZ$~is the magnetic energy. The $\DQZ$~frame (or $\dqZ$~frame for the IM) is chosen, because it yields simpler equations. We then invoke the two main simplifying assumptions made in basic modeling of AC motors:
\begin{itemize}
	\item the magnetic circuit is unsaturated, i.e., the flux-current relations are linear. This is equivalent to~$\HmDQZ$ (or $\HmdqZ$~for the IM) being at most quadratic with respect to the flux linkages
	\item the windings are sinusoidally distributed. This is equivalent to~$\HmDQZ$ (or $\HmdqZ$~for the IM) being independent of~$\te$. 
\end{itemize}
The most general magnetic energy corresponding to these assumptions is therefore
\ba\label{eqn:linear:general}
& & \HmDQZ(\phisDQZ, \phirDQZ) := a + \transpose{b}\!\!\phisDQZ + \transpose{c}\!\!\phirDQZ \ane \\
& & \qquad +\>\Transpose{\phisDQZ}\!\!D\phisDQZ + \Transpose{\phirDQZ}\!\!E\phisDQZ + \Transpose{\phirDQZ}\!\!F\phirDQZ \ans 
\ea
(or a similar expression of~$\HmdqZ$ for the IM), where 
$a$ is a constant, $b$ and~$c$ are $3 \times 1$~constant vectors and $D$, $E$ and~$F$ are $3 \times 3$~constant matrices;
without loss of generality, $D$ and $F$ (but not $E$) are assumed symmetric. Notice $a$~does not enter the equations~\eqnref{meca:hamiltonian} of analytical mechanics, hence can be freely chosen; the magnetic energy is thus parametrized by $3 + 3 + 6 + 9 + 6 = 27$~coefficients. By merely using the symmetries of \secref{symmetries} and the constraints of \secref{connections}, many of these $27$~coefficients are in fact linked, which will directly yield the classical unsaturated sinusoidal models.
Notice these models turn out to split into a $\DQ$- (or $\dq$-) part and a $0$-part, even without assuming a star connection, which is typical of unsaturated models.

\subsection{Synchronous Reluctant Motor}\label{sec:linear:srm}
As the SynRM has no rotor windings, its magnetic energy is independent of the rotor flux linkage according to~\secref{connection:norotor}; this implies that in~\eqnref{linear:general}, $c = \zero{3}{1}$ and $E = F = \zero{3}{3}$. On the other hand, the stator enjoys the symmetry of~\secref{symmetries:stator:rev}; by~\eqnref{symmetries:stator:rev:DQZ}, the coefficients in the odd powers of~$\phisZ$ are therefore zero, i.e., $\ind{b}{3}{} = \ind{D}{1}{3} = \ind{D}{2}{3} = \ind{D}{3}{1} = \ind{D}{3}{2} = 0$. Moreover, the rotor is symmetric with respect to two orthogonal planes containing the $D$\nobreakdash- and $Q$\nobreakdash-axes, hence by \eqnref{symmetries:rotor:swapQ:DQZ} and~\eqnref{symmetries:rotor:swapD:DQZ}, the coefficients in the odd powers of~$\phisD$ and~$\phisQ$ are zero, i.e., $\ind{b}{1}{} = \ind{b}{2}{} = \ind{D}{1}{2} = \ind{D}{2}{1} = 0$. 

Choosing $a:=0$, the simplest (pseudo)-Hamiltonian for the SynRM then reads in the $\DQZ$~frame
\beno
\HDQZ(\phisDQZ, \te, \ke) = \frac{\ke^2}{2\Jl\np^2}+\frac{1}{2}\GsD{\phisD}^2 + \frac{1}{2}\GsQ{\phisQ}^2 + \frac{1}{2}\GsZ{\phisZ}^2,
\eeno
where the coefficients have been renamed
\beno
\GsD := 2\ind{D}{1}{1} \qquad \GsQ := 2\ind{D}{2}{2} \qquad \GsZ := 2\ind{D}{3}{3}.
\eeno
Specializing \eqnref{DQZ:state-form}-\eqnref{DQZ:algebraic} to this case yields
\ba\label{eqn:linear:synrm:state-form}
\Tderive{\phisDQ} &=& \usDQ - \Rs\isDQ - \JJ2\we \phisDQ \ase \label{eqn:linear:synrm:state-form:DQ} \\
\Tderive{\phisZ} &=& \usZ - \Rs\isZ \ase \label{eqn:linear:synrm:state-form:Z}\\
\Tderive{\te} &=& \we \ase \\
\frac{\Jl}{\np} \Tderive{\we} &=& \np (\GsQ - \GsD)\phisD\phisQ - \Tl, \ase
\ea
with flux-current relations
\ba\label{eqn:linear:synrm:algebraic}
\isD &=& \GsD\phisD \ase \label{eqn:linear:synrm:algebraic:D}\\
\isQ &=& \GsQ\phisQ \ase \label{eqn:linear:synrm:algebraic:Q}\\
\isZ &=& \GsZ\phisZ. \ase \label{eqn:linear:synrm:algebraic:Z}
\ea
This is the classical model of the unsaturated sinusoidal SynRM. When the motor is star-connected, \eqnref{linear:synrm:state-form:Z} and \eqnref{linear:synrm:algebraic:Z} imply $0 = \isZ = \phisZ = \usZ$, which is used in the reduction from 3~to 2~axes of the classical model derivation of the model.

\subsection{Permanent Magnet Synchronous Motor}\label{sec:linear:pmsm}
Like the SynRM, the PMSM has no rotor windings, hence $c = \zero{3}{1}$ and $E = F = \zero{3}{3}$; similarly, \eqnref{symmetries:stator:rev:DQZ} imposes $\ind{b}{3}{} = \ind{D}{1}{3} = \ind{D}{2}{3} = \ind{D}{3}{1} = \ind{D}{3}{2} = 0$. By construction, the geometric saliency is aligned with the permanent magnet flux, which means the rotor is symmetric with respect to the plane defined by the $D$- and $0$- axes (provided the $D$-axis is chosen aligned with the permanent magnet flux); therefore, \eqnref{symmetries:rotor:swapQ:DQZ} implies $\ind{b}{2}{} = \ind{D}{1}{2} = \ind{D}{2}{1} = 0$.

Choosing $a:=\frac{\ind{b}{1}{}^2}{4\ind{D}{1}{1}}$, the simplest (pseudo)-Hamiltonian for the PMSM then reads in the $\DQZ$~frame
\ba\label{eqn:linear:pmsm}
& &\HDQZ(\phisDQZ, \te, \ke) = \frac{\ke^2}{2\Jl\np^2} \ane \\
& &\qquad+\>\frac{1}{2}\GsD(\phisD - \phiM)^2 + \frac{1}{2}\GsQ{\phisQ}^2 + \frac{1}{2}\GsZ{\phisZ}^2
\ea
where the coefficients have been renamed
\beno
\phiM = -\frac{\ind{b}{1}{}}{2\ind{D}{1}{1}}
\eeno
\beno
\GsD := 2\ind{D}{1}{1} \qquad \GsQ := 2\ind{D}{2}{2} \qquad \GsZ := 2\ind{D}{3}{3}.
\eeno
Notice $\phiM$~represents the magnitude of the permanent magnet flux.
Specializing \eqnref{DQZ:state-form}-\eqnref{DQZ:algebraic} to this case yields
\ba\label{eqn:linear:pmsm:state-form}
\Tderive{\phisDQ} &=& \usDQ - \Rs\isDQ - \JJ2\we \phisDQ \ase \label{eqn:linear:pmsm:state-form:DQ} \\
\Tderive{\phisZ} &=& \usZ - \Rs\isZ \ase \label{eqn:linear:pmsm:state-form:Z}\\
\Tderive{\te} &=& \we \ase \\
\frac{\Jl}{\np} \Tderive{\we} &=& \np\GsQ\phisQ\phiM + \np (\GsQ - \GsD)\phisD\phisQ - \Tl, \ans \ase
\ea
with flux-current relations
\ba\label{eqn:linear:pmsm:algebraic}
\isD &=& \GsD(\phisD - \phiM) \ase \label{eqn:linear:pmsm:algebraic:D}\\
\isQ &=& \GsQ\phisQ \ase \label{eqn:linear:pmsm:algebraic:Q}\\
\isZ &=& \GsZ\phisZ. \ase \label{eqn:linear:pmsm:algebraic:Z}
\ea
This is the classical model of the (salient) unsaturated sinusoidal PMSM; the non-salient PMSM is obtained with $\GsD = \GsQ = \Gs$, in which case the electro-magnetic torque reduces to $\Te = \np\GsQ\phisQ\phiM$. When the motor is star-connected, \eqnref{linear:pmsm:state-form:Z} and \eqnref{linear:pmsm:algebraic:Z} imply $0 = \isZ = \phisZ = \usZ$, which is used in the reduction from 3~to 2~axes of the classical model. 

Notice that when $\phiM = 0$, \eqnref{linear:pmsm:state-form}-\eqnref{linear:pmsm:algebraic} boil down the unsaturated sinusoidal SynRM~\eqnref{linear:synrm:state-form}-\eqnref{linear:synrm:algebraic}.

\subsection{Induction motor}\label{sec:linear:im}
The stator enjoys the symmetry of \secref{symmetries:stator:perm}; as $\HdqZ$~does not depend on~$\te$, \eqnref{symmetries:stator:perm:dqZ} reads
\beno
\HdqZ(\phisdqZ, \phirdqZ, \ke) = \HdqZ\bigl(\Rot3(\eta)\phisdqZ, \Rot3(\eta)\phirdqZ,\ke\bigr)
\eeno
for $\eta \in \{ \frac{2\pi}{3}, \frac{4\pi}{3}\}$,
which implies 
\bano
\transpose{b}\Rot3(\eta) &=& \transpose{b} \\
\transpose{c}\Rot3(\eta) &=& \transpose{c} \\ 
\transpose{\Rot3}(\eta)D\Rot3(\eta) &=& D \\
\transpose{\Rot3}(\eta)E\Rot3(\eta) &=& E \\
\transpose{\Rot3}(\eta)F\Rot3(\eta) &=& F.
\eano
As a consequence, $\ind{b}{1}{} = \ind{b}{2}{} = \ind{c}{1}{} = \ind{c}{2}{} = 0$, and $D$, $E$ and~$F$ have the form 
\beno
\begin{pmatrix}
	\xi & -\zeta & 0 \\
	\zeta & \xi & 0 \\
	0 & 0 & \chi
\end{pmatrix};
\eeno
as $D$~and $F$~are on the other hand symmetric, they are therefore diagonal. Moreover, $\ind{b}{3}{} = \ind{E}{3}{3} = 0$~by the symmetry condition~\eqnref{symmetries:stator:rev:dqZ}, and $\ind{c}{3}{} = 0$~by the symmetry condition~\eqnref{symmetries:rotor:rev:dqZ}; finally, the rotor is symmetric with respect to a plane (and even with respect to 3 planes, as we consider the rotor composed of 3 identical windings), hence $\ind{E}{2}{1} = \ind{E}{1}{2} = 0$ by the symmetry condition~\eqnref{symmetries:rotor:swapQ:DQZ}.

Choosing $a:=0$, the simplest (pseudo)-Hamiltonian for the IM then reads in the $\dqZ$~frame
\bano
& & \HdqZ(\phisdqZ, \phirdqZ, \te, \ke) = \frac{\ke^2}{2\Jl\np^2} + \frac{1}{2}\GlsZ{\phisZ}^2 + \frac{1}{2}\GlrZ{\phirZ}^2 \ane \\
& & \qquad+\>\frac{1}{2}\transpose{(\phisdq + \phirdq)}\Gm(\phisdq + \phirdq) \ane \\
& & \qquad+\>\frac{1}{2}\Transpose{\phisdq}\Gls\phisdq + \frac{1}{2}\Transpose{\phirdq}\Glr\phirdq,
\eano
where the coefficients have been renamed
\bano
\Gm &:=& \ind{E}{1}{1} = \ind{E}{2}{2} \\
\Gls &:=& 2\ind{D}{1}{1} - \ind{E}{1}{1} = 2\ind{D}{2}{2} - \ind{E}{2}{2} \\
\Glr &:=& 2\ind{F}{1}{1} - \ind{E}{1}{1} = 2\ind{F}{2}{2} - \ind{E}{2}{2} \\
\GlsZ &:=& 2\ind{D}{3}{3} \\
\GlrZ &:=& 2\ind{F}{3}{3}.
\eano
Specializing \eqnref{dqZ:state-form}-\eqnref{dqZ:algebraic} to this case yields
\ba\label{eqn:linear:im:state-form}
\Tderive{\phisdq} &=& \usdq - \Rs\isdq - \JJ2\ws \phisdq \ase \label{eqn:linear:im:state-form:stator:dq} \\
\Tderive{\phisZ} &=& \usZ - \Rs\isZ \ase \label{eqn:linear:im:state-form:stator:Z}\\
\Tderive{\phirdq} &=& -\Rr\irdq - \JJ2(\ws - \we) \phirdq \ase \label{eqn:linear:im:state-form:rotor:dq} \\
\Tderive{\phirZ} &=& - \Rr\irZ \ase \label{eqn:linear:im:state-form:rotor:Z}\\
\Tderive{\te} &=& \we \ase \\
\frac{\Jl}{\np} \Tderive{\we} &=& \np \Transpose{\phirdq}\JJ2\irdq - \Tl, \ase
\ea
with flux-current relations
\ba\label{eqn:linear:im:algebraic}
\isdq &=& \Gm(\phisdq + \phirdq) + \Gls\phisdq \ase \label{eqn:linear:im:algebraic:stator:dq} \\
\isZ &=& \GlsZ \phisZ \ase \label{eqn:linear:im:algebraic:stator:Z} \\
\irdq &=& \Gm(\phisdq + \phirdq) + \Glr\phirdq \ase \label{eqn:linear:im:algebraic:rotor:dq} \\
\irZ &=& \GlrZ \phirZ. \ase \label{eqn:linear:im:algebraic:rotor:Z}
\ea
This is the classical model of the unsaturated sinusoidal IM. When the motor is star-connected, \eqnref{linear:im:state-form:stator:Z} and \eqnref{linear:im:algebraic:stator:Z} imply $0 = \isZ = \phisZ = \usZ$; since the rotor windings are short-circuited, \eqnref{linear:im:state-form:rotor:Z} and \eqnref{linear:im:algebraic:rotor:Z} imply $0 = \irZ = \phirZ$; this is used in the reduction from 3 to 2~axes of the classical model.

%% file: experiments.tex
\section{Experimental validation on a PMSM}\label{sec:PMSM}
\begin{table}
	\setlength{\extrarowheight}{2pt}
	\caption{Specifications of the motor used in tests}\label{tbl:pmsm}
	\centering
	\begin{tabular}{lc}
		\firsthline
		Rated power & \SI{1500}{\watt}\\
		Rated speed & \SI{3000}{\rpm}\\
		Rated torque & \SI{6.06}{\newton\meter}\\
		Rated current (peak) & \SI{5.19}{\ampere}\\
		Rated voltage (peak) & \SI{245}{\volt}\\ \hline
		Pole pairs $\np$ & \num{5}\\
		Moment of inertia $\Jl$ & \SI{5.3e-3}{\kilogram.\meter^2}\\
		Stator resistance $\Rs$ & \SI{2.1}{\ohm}\\
		Magnet field $\phiM$ & \SI{0.155}{\weber}\\
		\lasthline
	\end{tabular}
\end{table}

We now validate the theory developed in the previous sections by experiments on an actual Surface-Mounted PMSM with Surface-Mounted Magnets (SPMSM). The motor (BMP1002F, see rated parameters in \tblref{pmsm}) is a high-end motor specifically designed for motion control, hence exhibit rather small non-sinusoidal and saturation effects; it was nevertheless chosen, because the star-point is accessible. Besides the motor, the experimental setup comprises a power stage taken from a Schneider Electric ATV71 \SI{1.5}{\kW} drive, a \SI{4}{\kW} load machine, and a dSpace environment. The motor is star-connected in all the experiments.

As in \secref{models}, we consider the simplest mechanical energy $\frac{\ke}{2\Jl\np^2}$ and a magnetic energy~$\HmDQZ$ independent of the kinetic momentum~$\ke$. Following \secref{connections:summary}, we can define a magnetic energy function $\HmzsDQ(\phisDQ, \te)$ independent of the $0$-variables, yielding the unconstrained $\DQ$-subsystem 
	\ba\label{eqn:nonsin:state}
	\Tderive{\phisDQ} &=& \vsDQ - \Rs\isDQ - \JJ2\we\phisDQ \ase\\
	\Tderive{\te} &=& \we \ase\\
	\frac{\Jl}{\np} \Tderive{\we} &=& \Te - \Tl, \ase
	\ea
with constitutive relations
	\ba\label{eqn:nonsin:constit}
	\isDQ &=& \Pderive{\HmzsDQ}{\phisDQ}(\phisDQ, \te) \ase\\
	\Te &=& -\np\Pderive{\HmzsDQ}{\te}(\phisDQ, \te) + \np\Transpose{\isDQ}\JJ2\phisDQ.\ase
	\ea
To study non-sinusoidal effects, we will also consider the constrained $0$-stator subsystem, written as the purely differential version of~\eqnref{connection:star:zero} in the form
	\ba\label{eqn:nonsin:phiz}
	\Tderive{\phisZ} &=& \vsZ - v_\starpoint \ase \\
	&=& - (\Jac{}{\phisZ} \Pderive{\HmzDQZ}{\phisZ} )^{-1}f(\phisDQ, \phisZ, \te; \vsxy, \Tl); \ans \ase
	\ea
it involves the full magnetic energy~$\HmzDQZ(\phisDQZ,\te)$.

\subsection{Non-sinusoidal model}\label{sec:nonsin}
By the stator symmetry~\eqnref{symmetries:stator:perm:DQZ}, $\HmzDQZ$~is $\frac{2\pi}{3}$\nobreakdash-periodic with respect to $\te$, hence can be expanded as the Fourier series
\bano
& & \HmzDQZ(\phisDQZ, \te) = \HbmZDQZ(\phisDQZ) \\
& & \qquad+\>\sum_{k=1}^{\infty} \fcos{3k}(\phisDQZ)\cos{3k\te} + \fsin{3k}(\phisDQZ)\sin{3k\te}.
\eano
This implies that $v_\starpoint$, $\isDQZ$ and $\Te$, which are derived from $\HmzDQZ$, may contain only harmonics of order $3k$.
Besides, the stator symmetry~\eqnref{symmetries:stator:rev:DQZ} imposes
\bano
\HbmZDQZ(\phisDQ, \phisZ) &=& \HbmZDQZ(\phisDQ, -\phisZ) \\
\fcos{3k}(\phisDQ, \phisZ) &=& (-1)^k\fcos{3k}(\phisDQ, -\phisZ) \\
\fsin{3k}(\phisDQ, \phisZ) &=& (-1)^k\fsin{3k}(\phisDQ, -\phisZ).
\eano

Similarly, $\HmzsDQ$~can be expanded as 
\ba\label{eqn:nonsin:fourrier:star}
& & \HmzsDQ(\phisDQ, \te) = \HbmZsDQ(\phisDQ) \ane \\
& & \qquad+\>\sum_{k=1}^{\infty} \fscos{3k}(\phisDQ)\cos{3k\te} + \fssin{3k}(\phisDQ)\sin{3k\te},
\ea
and
\bano
\HbmZsDQ(\phisDQ) &=& \HbmZsDQ(\phisDQ)\\
\fscos{3k}(\phisDQ) &=& (-1)^k\fscos{3k}\\
\fssin{3k}(\phisDQ) &=& (-1)^k\fssin{3k}(\phisDQ).
\eano
The last two relations mean $\fscos{6k + 3}=\fssin{6k + 3}=0$,
implying that $\isDQZ$ and $\Te$, which are derived from $\HmzsDQ$, may contain in fact only harmonics of order $6k$.
Torque ripple at $6^{th}$ times the velocity is indeed a well-known phenomenon for AC motors, see e.g.~\cite[chap.~2]{Bose2002book}.

\begin{figure}
	\centering
	\subfloat[Time evolution of~$\Te$.\label{fig:nonsin:te:time}]{\includegraphics[width=8cm]{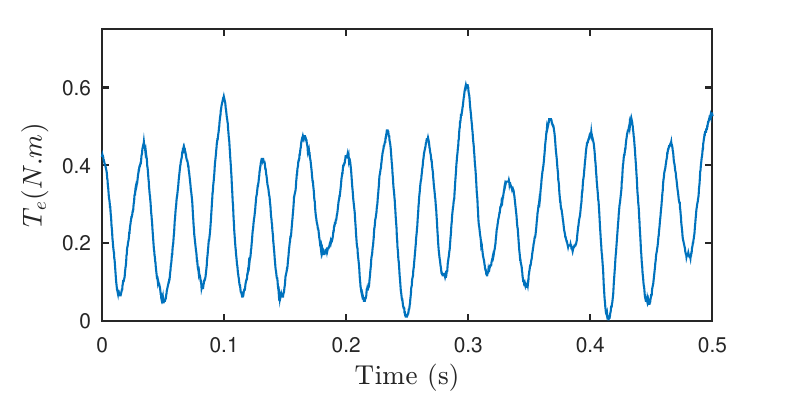}} \\
	\subfloat[Spectrum of~$\Te$.\label{fig:nonsin:te:freq}]{\includegraphics[width=8cm]{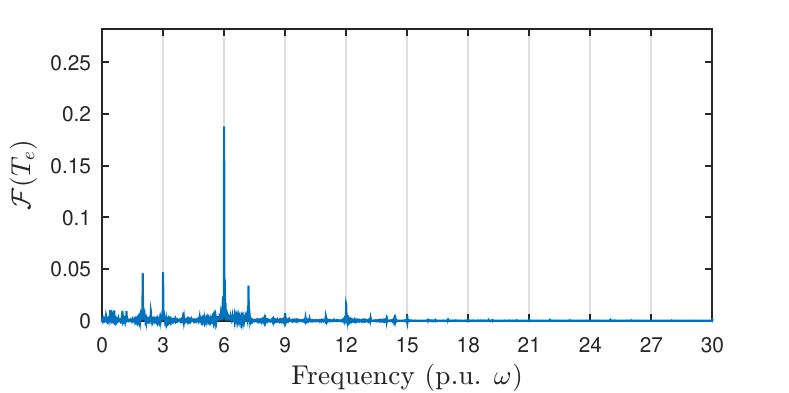}}
	\caption{Electromagnetic torque $\Te$.}
	\label{fig:nonsin:te}
\end{figure}

\begin{figure}
	\centering
	\subfloat[Time evolution of~$v_\starpoint$.\label{fig:nonsin:vn:time}]{\includegraphics[width=8cm]{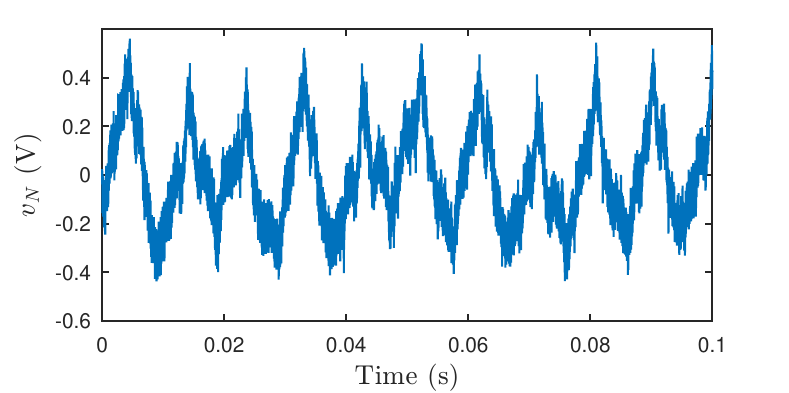}} \\
	\subfloat[Spectrum of~$v_\starpoint$.\label{fig:nonsin:vn:freq}]{\includegraphics[width=8cm]{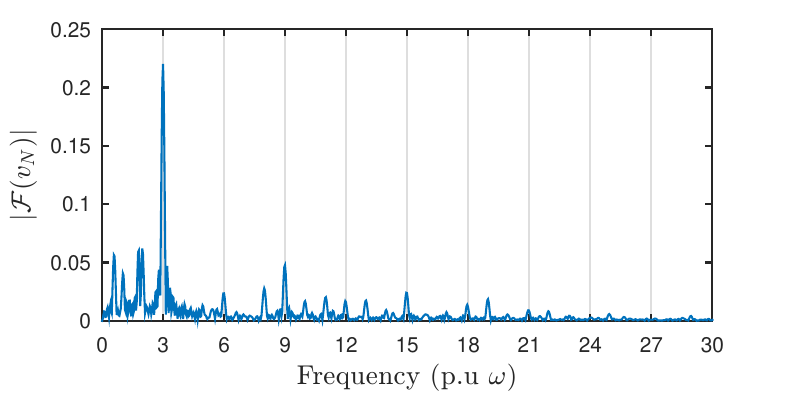}}
	\caption{Potential $v_\starpoint$ of the star point.}
	\label{fig:nonsin:vn}
\end{figure}

We experimentally checked the previous conclusions about~$\Te$ and~$v_\starpoint$:
\begin{itemize}
	\item in a first experiment, the motor was excited with a \SI{5}{\hertz} sinusoidal voltage with the load machine disconnected. The time evolution of~$\Te$ and its spectrum are displayed in \figref{nonsin:te}. The main harmonic at $6\we$ is clearly visible; the small-amplitude harmonics at~$2\we$ and~$3\we$, not explained by our model, are probably due to other imperfections of the motor or the power stage
	\item in a second experiment, the motor (still star-connected) was disconnected from the drive, and driven at the (electrical) speed~$\we:=\SI{35}{\hertz}$ by the load machine. The phase voltage drops $\vsabc - v_\starpoint = \Tderive{\phisabc}$ were recorded, and $v_\starpoint$ was computed according to~\eqnref{nonsin:phiz} as $v_\starpoint = \frac{\vsZ}{\sqrt{3}} - \frac{1}{3}\bigl(\Tderive{\phisa} + \Tderive{\phisb} + \Tderive{\phisc}\bigr)$, with $\vsZ$ an arbitrary constant ($0$ in our case); being disconnected from the drive allows us to recover the rather small $v_\starpoint$, which would otherwise be partly hidden by the imperfections of the power stage. The time evolution of $v_\starpoint$ and its spectrum  are displayed in \figref{nonsin:vn}. As anticipated, the main harmonic is at $3\we$.
\end{itemize}

\subsection{Saturated model}\label{sec:nonlin}
The ripple caused by non-sinusoidal effects is usually small with respect to the rated values, and can be neglected in the model for control proposes. We therefore concentrate on the fundamental term~$\HbmZsDQ$ in the Fourrier expansion of the magnetic energy~\eqnref{nonsin:fourrier:star};
as in the unsaturated case, the rotor is symmetric with respect to the plane defined by the $D$-and $0$\nobreakdash-axes, so $\HbmZsDQ$~is even with respect to~$\phisQ$ by \eqnref{symmetries:rotor:swapQ:DQZ}. Following~\cite{JebaiMMR2016IJoC}, we consider a perturbation of the unsaturated model~\eqnref{linear:pmsm} in the form
\ba\label{eqn:nonlinear:pmsm}
& & \HbmZsDQ(\phisD, \phisQ) = \frac{1}{2}\fD(\phisD - \phiM) + \frac{1}{2}\fQ\bigl(\phisQ^2\bigr) \ane \\
& & \qquad+\>\frac{1}{2}\fx\bigl(\phisD - \phiM, \phisQ^2\bigr),
\ea
where $\fD$, $\fQ$, $\fx$ are polynomials of order at most 4,
\bano
\fD(\psisD) &:=& \GsD\Bigl(\psisD^2 + \frac{\psisD^3}{6\phi_1^D} + \frac{\psisD^4}{12{\phi_2^D}^2}\Bigr) \\
\fQ(\phisQ^2) &:=& \GsQ\Bigl(\phisQ^2 + \frac{\phisQ^4}{12{\phi_1^Q}^2}\Bigr) \\
\fx(\psisD, \phisQ^2) &:=& \GsD\Bigl(\frac{\psisD}{2\phi_1^X} + \frac{\psisD^2}{{\phi_2^X}^2} \Bigr)\phisQ^2,
\eano
and~$\psisD := \phisD - \phiM$. The 7 parameters $\GsD$, $\GsQ$, $\phi_1^D$, $\phi_2^D$, $\phi_1^Q$, $\phi_1^X$ and $\phi_2^X$ have to be determined experimentally. This model can be seen as a fourth-order Taylor expansion of the true energy function, whereas the unsaturated model is just a second-order expansion.

The 7 parameters were experimentally identified with signal injection, see \cite{JebaiMMR2016IJoC} for details about the procedure; the values are given in \tblref{nonlin:pmsm:params}. \figref{nonlin:pmsm:current-flux} 
illustrate the good agreement between the model and the experimental data; the saturation effects are clearly visible, though this high-end motor exhibits little saturation compared to more standard motors.

\subsection{The role of magnetic saturation for sensorless control} \label{sec:nonlin:hf}
In sensorless control, the actual control input in the model~\eqnref{nonsin:state}-\eqnref{nonsin:constit} is $\vsabc$, or equivalently $\vsAB$, and the actual measurement is $\isabc$, or equivalently $\isAB$. At zero velocity, there is a loss of observability, which makes the control problem more difficult in this zone. A fairly recent but now well-established method, called signal injection, to overcome this problem is to use an hybrid control law of the form
\bano
\vsAB &:=&  \lf{\vsAB} + \hf{\vsAB}f(t),
\eano
where $\lf{\vsAB}$ and $\hf{\vsAB}$ are ``slowly-varying'' inputs, and $f(t)$ is a ``high-frequency'' periodic signal. Following the analysis based on second-order averaging of~\cite{CombeJMMR2016ACC,JebaiMMR2016IJoC}, the net effect of this signal injection is to make available for use in the control law the ``virtual output''
\bano
\hf{\isAB} &:=& \Jac{}{\phisDQ}\isAB\cdot\Rot2(-\te)\hf{\vsAB}\\
&=&\Jac{}{\phisDQ}\bigl(\Rot2(\te)\Pderive{\HDQ}{\phisDQ}(\phisDQ)\bigr)\cdot\Rot2(-\te)\hf{\vsAB}\\
&=&\Rot2(\te)\Jac{}{\phisDQ}\Pderive{\HDQ}{\phisDQ}(\phisDQ)\Rot2(-\te)\hf{\vsAB},
\eano
where
\bano
S(\te,\phisDQ) &:=& \Rot2(\te)\Jac{}{\phisDQ}\Pderive{\HDQ}{\phisDQ}(\phisDQ)\Rot2(-\te)
\eano
is the so-called saliency matrix; thanks to this extra output, the observability problem disappears provided $S$ effectively depends on~$\te$.

For an unsaturated PMSM described by~\eqnref{linear:pmsm}, the saliency matrix boils down to
\bano
&&\Rot2(\te)\begin{pmatrix}\GsD&0\\ 0&\GsQ\end{pmatrix}\Rot2(-\te).
\eano
If the motor has little geometric saliency, as is the case for an SPMSM, $\GsD\approx\GsQ$; this implies $S\approx\GsD I$, i.e., $S$ does not depend on~$\te$. Ignoring magnetic saturation would lead to the wrong conclusion that such a motor cannot be controlled by signal injection. But when saturation is taken into account through~\eqnref{nonlinear:pmsm}, $S$ will in general depend on~$\te$ thanks in particular to the off-diagonal terms of $\Jac{}{\phisDQ}\Pderive{\HDQ}{\phisDQ}$. Correctly describing saturation is thus paramount for operation at low velocity, which is experimentally confirmed, see e.g.~\cite{JebaiMMR2016IJoC}.


%
%
\begin{figure}
	\centering
	\setlength\figurewidth{200pt}%
	\setlength\figureheight{150pt}%
	\input{Figures/spm_current_flux.tikz}%
 \\
	\caption{Experimental current-flux relations~$\isD(\phisD - \phiM, 0)$~(\textcolor{blue}{---}) and~$\isQ(0, \phisQ)$~(\textcolor{green!50!black}{---}) vs model (\textcolor{red}{--{ }--}).}
	\label{fig:nonlin:pmsm:current-flux}
\end{figure}
%

\begin{table}
	\setlength{\extrarowheight}{2pt}
	\caption{Estimated parameters}
	\centering
	\begin{tabular}{lclclc}
		\firsthline
		\multicolumn{2}{l}{Parameters of $\fD$} & \multicolumn{2}{l}{Parameters of $\fQ$} & \multicolumn{2}{l}{Parameters of $\fx$}\\ \hline
		$\frac{1}{\GsD}$ & \SI{8.8}{\milli\henry} & $\frac{1}{\GsQ}$ & \SI{7.7}{\milli\henry} && \\
		$\phi_1^D$ & \SI{0.533}{\weber} & $\phi_1^Q$ & \SI{0.228}{\weber} & $\phi_1^x$ & \SI{0.116}{\weber} \\
		$\phi_2^D$ & \SI{0.200}{\weber} &&& $\phi_2^x$ & \SI{0.111}{\weber} \\
		\lasthline
	\end{tabular}
	\label{tbl:nonlin:pmsm:params}
\end{table}